# WISHART DISTRIBUTIONS FOR DECOMPOSABLE GRAPHS

By Gérard Letac and Hélène Massam[1]

*Université Paul Sabatier and York University*

When considering a graphical Gaussian model $\mathcal{N}_G$ Markov with respect to a decomposable graph $G$, the parameter space of interest for the precision parameter is the cone $P_G$ of positive definite matrices with fixed zeros corresponding to the missing edges of $G$. The parameter space for the scale parameter of $\mathcal{N}_G$ is the cone $Q_G$, dual to $P_G$, of incomplete matrices with submatrices corresponding to the cliques of $G$ being positive definite. In this paper we construct on the cones $Q_G$ and $P_G$ two families of Wishart distributions, namely the Type I and Type II Wisharts. They can be viewed as generalizations of the hyper Wishart and the inverse of the hyper inverse Wishart as defined by Dawid and Lauritzen [*Ann. Statist.* **21** (1993) 1272–1317]. We show that the Type I and II Wisharts have properties similar to those of the hyper and hyper inverse Wishart. Indeed, the inverse of the Type II Wishart forms a conjugate family of priors for the covariance parameter of the graphical Gaussian model and is strong directed hyper Markov for every direction given to the graph by a perfect order of its cliques, while the Type I Wishart is weak hyper Markov. Moreover, the inverse Type II Wishart as a conjugate family presents the advantage of having a multidimensional shape parameter, thus offering flexibility for the choice of a prior.

Both Type I and II Wishart distributions depend on multivariate shape parameters. A shape parameter is acceptable if and only if it satisfies a certain eigenvalue property. We show that the sets of acceptable shape parameters for a noncomplete $G$ have dimension equal to at least one plus the number of cliques in $G$. These families, as conjugate families, are richer than the traditional Diaconis–Ylvisaker conjugate families which all have a shape parameter set of dimension one. A decomposable graph which does not contain a three-link chain as an induced subgraph is said to be homogeneous. In this case, our Wisharts are particular cases of the Wisharts on homogeneous cones

Received January 2005; revised July 2006.
[1]Supported by NSERC Grant A8947.

as defined by Andersson and Wojnar [*J. Theoret. Probab.* **17** (2004) 781–818] and the dimension of the shape parameter set is even larger than in the nonhomogeneous case: it is indeed equal to the number of cliques plus the number of distinct minimal separators. Using the model where $G$ is a three-link chain, we show by computing a 7-tuple integral that in general we cannot expect the shape parameter sets to have dimension larger than the number of cliques plus one.

**1. Introduction.** The primary aim of this paper is to develop a new family of conjugate prior distributions with attractive Markov properties for the covariance parameter, or equivalently the precision parameter, of graphical Gaussian models Markov with respect to a decomposable graph $G$. While doing so, we are led to define two new classes of Wishart distributions and their inverses and to study their properties.

Let us recall that an undirected graph is a pair $(V, \mathcal{E})$ where $V = \{1, \ldots, r\}$ and $\mathcal{E}$ is a family of subsets $\{i, j\}$ of $V$ of size 2. It will be convenient to consider the set $E \subset V \times V$ of $(i, j)$ such that either $i = j$ or $\{i, j\}$ is in $\mathcal{E}$, rather than $\mathcal{E}$ and, since $E$ and $\mathcal{E}$ carry the same information, to speak about the graph $G = (V, E)$. Any $(i, j)$ such that $i \neq j$ will be called an edge. An $r$-dimensional Gaussian model is said to be Markov with respect to $G$ if for any edge $(i, j)$ not in $E$, the $i$th and $j$th variables are conditionally independent given all the other variables. Such models are known as covariance selection models (see [8]) or graphical Gaussian models (see [18] or [11]). Without loss of generality, we can assume that these models are centered $N_r(0, \Sigma)$, and it is well known that they are characterized by the parameter set $P_G$ of the precision matrices, which is the set of positive definite matrices $K = \Sigma^{-1}$ such that $K_{ij} = 0$ whenever the edge $(i, j)$ is not in $E$. Equivalently, if we denote by $M$ the linear space of symmetric matrices of order $r$, by $M_r^+ \subset M$ the cone of positive definite (abbreviated $> 0$) matrices, by $I_G$ the linear space of symmetric incomplete matrices $x$ with missing entries $x_{ij}, (i, j) \notin E$, and by $\pi: M \mapsto I_G$ the projection of $M$ into $I_G$, the parameter set of the Gaussian model can be described as the set of incomplete matrices $\Sigma_G = \pi(\Sigma)$ with $\Sigma = K^{-1}$ and $K \in P_G$. Indeed it is easy to verify that the entries $\Sigma_{ij}, (i, j) \notin E$ are such that

$$\Sigma_{ij} = \Sigma_{i, V \setminus \{i, j\}} \Sigma_{V \setminus \{i, j\}, V \setminus \{i, j\}}^{-1} \Sigma_{V \setminus \{i, j\}, j},$$

and are therefore not free parameters of the Gaussian models. One can prove that the correspondence between $K$ and the incomplete matrix $\Sigma_G = \pi(\Sigma)$ is one to one. We write $\Sigma_G = \varphi(K) = \pi(K^{-1})$. We note that $\varphi$ is not explicit when $G$ is not decomposable.

Henceforth in this paper, we will assume that $G$ is decomposable. The reader is referred to [11] for all the common notions of graphical models used in this paper. We will now simply recall some basic facts and traditional notation we will use throughout this paper. Every decomposable graph admits a perfect order of its cliques. Let $(C_1, \ldots, C_k)$ be such an order. We use the notation $H_1 = R_1 = C_1$, while for $j = 2, \ldots, k$ we write

$$H_j = C_1 \cup \cdots \cup C_j, \qquad R_j = C_j \setminus H_{j-1}, \qquad S_j = H_{j-1} \cap C_j.$$

The $S_j, j = 2, \ldots, k$, are the minimal separators of $G$. Some of these separators can be identical. We let $k' \leq k-1$ denote the number of distinct separators and $\nu(S)$ denote the multiplicity of $S$ that is the number of $j$ such that $S_j = S$. Lauritzen [11] has proven that the multiplicity $\nu(S)$ of a given minimal separator $S$ is positive and independent of the perfect order of the cliques considered.

For $G$ given decomposable with the set of cliques $\{C_1, \ldots, C_k\}$ and $\Sigma^{-1} \in P_G$, the incomplete matrix $\Sigma_G$ is completely determined by its submatrices $\{\Sigma_{C_i}, i = 1, \ldots, k\}$ where, of course, for each $i = 1, \ldots, k, \Sigma_{C_i}$ is positive definite. When considering the parameter space of the graphical Gaussian model corresponding to $G$ decomposable, we are therefore led to consider the two cones

(1.1) $$P_G = \{y \in M_r^+ | y_{ij} = 0, (i,j) \notin E\},$$

(1.2) $$Q_G = \{x \in I_G | x_{C_i} > 0, i = 1, \ldots, k\}.$$

Dawid and Lauritzen ([7], Section 7) defined two distributions on $Q_G$, namely, the hyper Wishart distribution as the distribution of the maximum likelihood estimator of $\Sigma_G$, and the hyper inverse Wishart distribution as the Diaconis–Ylvisaker conjugate prior distribution for $\Sigma_G$. Subsequently Roverato [16] gave the distribution of $K = \Sigma^{-1} = \varphi^{-1}(\Sigma_G)$ when $\Sigma_G$ follows the hyper inverse Wishart distribution. We will call this distribution of $K$ on $P_G$ the $G$-Wishart. The search for a rich and flexible class of conjugate prior distributions for $\Sigma_G$, or equivalently for $K = \Sigma^{-1}$, remains a topic of high interest to statisticians.

When $G$ is complete, $P_G = Q_G = M_r^+$ and we define the regular Wishart distribution on the cone of positive definite matrices of dimension $r = |V|$ by

$$\frac{1}{2^{rp}\Gamma_r(p)|\Sigma|^p} e^{-(1/2)\mathrm{tr}(x\Sigma^{-1})} |x|^p |x|^{-(r+1)/2} \mathbf{1}_{M_r^+}(x) \, dx,$$

where $p > \frac{r-1}{2}$ is the one-dimensional shape parameter and $\Sigma \in M_r^+$ is the scale parameter.

When $G$ is decomposable, the hyper and hyper inverse Wisharts have been constructed as a Markov combination (with respect to $G$) of the Wishart and its inverse, respectively, and so, like the Wishart, they have a one-dimensional shape parameter and a scale parameter in $Q_G$. Dawid and Lauritzen [7] have shown that these distributions have Markov properties: the hyper Wishart is weak hyper Markov while the hyper inverse Wishart is strong hyper Markov.

In this paper, we will construct a family of distributions, called Type I Wisharts, defined on $Q_G$, and another family, called Type II Wisharts, defined on $P_G$. We shall see in Section 4 that the inverses of the Type II distributions, like the hyper inverse Wisharts, form a family of conjugate prior distributions for the scale parameter of the graphical Gaussian model. We will also show that they are strong directed hyper Markov in the direction given to the graph $G$ by any choice of a perfect numbering of its vertices. This property is parallel to the strong hyper Markov property of the hyper inverse Wishart. We will also show that the Type I Wishart is weak hyper Markov, a property parallel to the weak hyper Markov property of the hyper Wishart. The attractive feature of the inverse Type II Wishart family of conjugate distributions is that, except in the trivial case where $G$ is complete, the set of shape parameters is of dimension strictly greater than the number $k$ of cliques in $G$, thus offering a flexible class of conjugate prior distributions for $\Sigma_G$. We shall also note in Section 4 that it forms a class of enriched conjugate priors for $\Sigma_G$ in the sense of [5].

To construct these two families, we define two natural exponential families of distributions affiliated with the Wishart, one on $Q_G$ and one on $P_G$. Let $(C_1, \ldots, C_k)$ denote a perfect order of the cliques of $G$ and let $(S_2, \ldots, S_k)$ be its corresponding sequence of minimal separators, some of them being possibly identical. We consider functions of the type

$$H_G(\alpha, \beta; x) = \frac{\prod_{i=1}^{k} |x_{C_i}|^{\alpha_i}}{\prod_{i=2}^{k} |x_{S_i}|^{\beta_i}}, \qquad x \in Q_G,$$

where $\alpha$ and $\beta$ are two real-valued functions on the collections $\mathcal{C}$ and $\mathcal{S}$ of cliques and separators, respectively, such that $\alpha(C_i) = \alpha_i, \beta(S_j) = \beta_j$ with $\beta_i = \beta_j$ if $S_i = S_j$. These functions play a very special role in the definition of the two families of distributions that we define. Indeed, if we let $c_i = |C_i|$ and $s_i = |S_i|$ denote the cardinality of $C_i$ and $S_i$, respectively, and if we denote

$$\mu_G(dx) = \frac{\prod_{i=1}^{k} |x_{C_i}|^{-(c_i+1)/2}}{\prod_{i=2}^{k} |x_{S_i}|^{-(s_i+1)/2}} \mathbf{1}_{Q_G}(x)\, dx,$$

the family of distributions we define on $Q_G$ is, for a given $(\alpha, \beta)$, the natural exponential family generated by

(1.3) $\quad H_G(\alpha, \beta; x)\mu_G(dx) = H_G(\alpha - \tfrac{1}{2}(c+1), \beta - \tfrac{1}{2}(s+1); x)\mathbf{1}_{Q_G}(x)\, dx.$

The measure (1.3) can be seen as a Markov combination generalization of the measure $|x|^p |x|^{-(r+1)/2} \mathbf{1}_{M_r^+}(x)\, dx$ generating the Wishart distribution, for a given $p > \frac{r-1}{2}$.

In Section 3 we will introduce the set $\mathcal{A}$ of $(\alpha, \beta)$ such that the following integral converges and satisfies

(1.4) $$\int_{Q_G} e^{-\mathrm{tr}(xy)} H_G(\alpha, \beta; x) \mu_G(dx) = \Gamma_\mathrm{I}(\alpha, \beta) H_G(\alpha, \beta; \varphi(y)),$$

where $\Gamma_\mathrm{I}(\alpha, \beta)$ is some function of $(\alpha, \beta)$ independent of $y \in P_G$. When $(\alpha, \beta)$ is in $\mathcal{A}$ we say that $H_G(\alpha, \beta; x)$ has the eigenvalue property with corresponding eigenvalue $\Gamma_\mathrm{I}(\alpha, \beta)$ and we define the Type I Wishart distribution on $Q_G$ as the distribution with density

$$\frac{1}{\Gamma_\mathrm{I}(\alpha, \beta) H_G(\alpha, \beta; \varphi(y))} e^{-\mathrm{tr}(xy)} H_G(\alpha, \beta; x) \mu_G(dx)$$

and with parameters $(\alpha, \beta, y)$. In a parallel way, we define a set $\mathcal{B}$ of $(\alpha, \beta)$ for which an eigenvalue property similar to (1.4) holds for the Type II Wishart distribution defined on $P_G$.

In order to fully describe the Type I and II Wishart distributions, it is then necessary to know the sets $\mathcal{A}$ and $\mathcal{B}$. In Section 3.2 we show that, for any $G$, the hyper Wishart and the $G$-Wishart are particular cases of Type I and II distributions, respectively. More precisely, we describe the sets $\mathcal{A}_1 \subset \mathcal{A}$ and $\mathcal{B}_1 \subset \mathcal{B}$ such that for $(\alpha, \beta) \in \mathcal{A}_1$, the Type I Wishart is the hyper Wishart and for $(\alpha, \beta) \in \mathcal{B}_1$, the Type II Wishart is the $G$-Wishart. In Section 3.3 we consider the particular class of decomposable graphs $G$ which do not contain the three-link chain, which we call $A_4$, as an induced subgraph. Such graphs are called homogeneous. When $G$ is homogeneous, we describe the sets $\mathcal{A}$ and $\mathcal{B}$ completely and show that they are open sets of dimension $k + k'$, the number of cliques plus the number of distinct separators in $G$. For $G$ homogeneous, the cones $Q_G$ and $P_G$ are homogeneous and we see that the Type I and II Wisharts then belong to the class of Wisharts on homogeneous cones defined by Andersson and Wojnar [3]. In Section 3.4 we consider nonhomogeneous graphs. In that case, we have, so far, only partial knowledge of $\mathcal{A}$ and $\mathcal{B}$. For each perfect order $P$ of the cliques, we define a $(k+1)$-dimensional subset $A_P$ of $\mathcal{A}$ such that, for $(\alpha, \beta) \in A_P$, (1.4) holds. We therefore know the subset $\bigcup A_P$ of $\mathcal{A}$, but not all of $\mathcal{A}$. Similarly we define a $(k+1)$-dimensional subset $B_P$ of $\mathcal{B}$ such that we know the subset $\bigcup B_P$ of $\mathcal{B}$, but not all of $\mathcal{B}$. We conjecture that the equalities $\mathcal{A} = \bigcup A_P$ and $\mathcal{B} = \bigcup B_P$ hold in general for nonhomogeneous graphs and that, thus, the dimension of the manifolds $\mathcal{A}$ and $\mathcal{B}$ is generally $k + 1 < k + k'$. In Section 3.4 we verify that these two equalities hold when $G = A_4$ by computing, in this case, the 7-tuple integral corresponding to (1.4).

In Section 4 we give the conjugacy and hyper Markov properties mentioned above. We also give the Laplace transforms of the Type I and II Wisharts and the expected values of the Type I, Type II and inverse Type II Wisharts. The necessary preliminaries for understanding the cones $P_G$ and $Q_G$ and the measures we define on them are given in Section 2.1. In Section 2.2 we give the results needed to work with homogeneous graphs. Most proofs are deferred to the Appendix.

## 2. Preliminaries.

2.1. *Measures on $P_G$ and $Q_G$.* For the graph $G = (V, E)$, $V = \{1, \ldots, r\}$, we write $i \sim j$ to indicate that the edge $\{i, j\}$ is in $\mathcal{E}$. An undirected graph $G$ is said to be decomposable if it does not contain a cycle of length greater than or equal to four as an induced subgraph and if it is connected. For all the notions related to decomposable graphs that we will introduce below, the reader is referred to [11], Chapter 2. We denote by $Z_G$ the real linear space of symmetric matrices $y$ of order $r$ such that $y_{ij} = 0$ if $(i, j) \notin E$. We denote by $I_G$ the real linear space of functions $(i, j) \mapsto x_{ij}$ from $E$ to $\mathbb{R}$ such that $x_{ij} = x_{ji}$. The elements of $I_G$ are called *G-incomplete symmetric matrices*. For a decomposable graph, we have defined the cones $P_G \subset Z_G$ and $Q_G \subset I_G$ in (1.1) and (1.2). Recall that $M_r^+$ denotes the cone of positive definite symmetric matrices of order $r$. Gröne et al. [10] proved the following.

PROPOSITION 2.1. *When $G$ is decomposable, for any $x$ in $Q_G$ there exists a unique $\hat{x}$ in $M_r^+$ such that for all $(i, j)$ in $E$ we have $x_{ij} = \hat{x}_{ij}$ and such that $\hat{x}^{-1}$ is in $P_G$.*

This defines a bijection between $P_G$ and $Q_G$,

$$\varphi : y = (\widehat{x})^{-1} \in P_G \mapsto x = \varphi(y) = \pi(y^{-1}) \in Q_G,$$

where $\pi$ denotes the projection of $M$ onto $I_G$. The explicit expression of $\widehat{x}^{-1}$ is given in (2.3) below. For $(x, y) \in I_G \times Z_G$, we write $\operatorname{tr}(xy) = \langle x, y \rangle = \sum_{(i,j) \in E} x_{ij} y_{ij}$. By Proposition 2.1, we have for $x \in Q_G$ $\langle x, y \rangle = \operatorname{tr}(\hat{x}y)$, where $\operatorname{tr}(\hat{x}y)$ is defined in the classical way. Thus although $xy$ does not make sense, the notation $\operatorname{tr}(xy)$ is quite convenient. We also use the following notation: if $C$ is a complete subset of vertices and if $x_C = (x_{ij})_{i,j \in C}$ is a matrix, we denote by $(x_C)^0 = (x_{ij})_{i,j \in V}$ the matrix such that $x_{ij} = 0$ for $(i, j) \notin C \times C$.

The following theorem gathers some basic results on decomposable graphs. Part 1 is due to Andersson [2], parts 2 and 3 can be found in [11], Chapter 5 and part 4 is due to Roverato [16].

THEOREM 2.1. *Let $G$ be a decomposable graph. Then:*

1. *The convex open cones $P_G$ and $Q_G$ are dual to each other in the sense that*

(2.1) $$P_G = \{y \in Z_G; \operatorname{tr}(xy) > 0 \ \forall x \in \overline{Q}_G \setminus \{0\}\},$$

(2.2) $$Q_G = \{x \in I_G; \operatorname{tr}(xy) > 0 \ \forall y \in \overline{P}_G \setminus \{0\}\}.$$

2. *For $x \in Q_G$ we have that $y = \hat{x}^{-1}$ is in $P_G$ and*

$$y = \sum_{C \in \mathcal{C}} (x_C^{-1})^0 - \sum_{S \in \mathcal{S}} \nu(S)(x_S^{-1})^0. \tag{2.3}$$

3. *For $x \in Q_G$ we have*

$$\det \hat{x} = \frac{\prod_{C \in \mathcal{C}} (\det x_C)}{\prod_{S \in \mathcal{S}} (\det x_S)^{\nu(S)}}.$$

4. *The absolute value of the Jacobian of the bijection $x \mapsto y = \hat{x}^{-1}$ from $Q_G$ to $P_G$ is*

$$\prod_{C \in \mathcal{C}} (\det x_C)^{-|C|-1} \prod_{S \in \mathcal{S}} (\det x_S)^{(|S|+1)\nu(S)}. \tag{2.4}$$

The proof of part 1 is given in the Appendix. For $G$ complete, part 4 above becomes the following.

LEMMA 2.1 ([14]). *The Jacobian of the change of variable $x \in M_r^+ \mapsto y = x^{-1} \in M_r^+$ is $|y|^{-(r+1)}$.*

We now introduce the measures which will be the generating measures of the new Wishart exponential families on $P_G$ and $Q_G$ that we are going to define in the next section. Let $\alpha$ and $\beta$ be two real valued functions on $\mathcal{C}$ and $\mathcal{S}$, respectively. An example of such functions $\alpha$ and $\beta$ is

$$C \in \mathcal{C} \mapsto \alpha(C) = |C| \quad \text{and} \quad S \in \mathcal{S} \mapsto \beta(S) = |S|.$$

We denote these examples $\alpha = c$ and $\beta = s$. Another example is, for a constant $p$ given,

$$C \in \mathcal{C} \mapsto \alpha(C) = p \quad \text{and} \quad S \in \mathcal{S} \mapsto \beta(S) = p,$$

simply denoted $\alpha = p$ and $\beta = p$. For $x \in Q_G$ we adopt the notation

$$H_G(\alpha, \beta; x) = \frac{\prod_{C \in \mathcal{C}} (\det x_C)^{\alpha(C)}}{\prod_{S \in \mathcal{S}} (\det x_S)^{\nu(S)\beta(S)}}. \tag{2.5}$$

The functions $H_G$ for the particular case $\alpha = -\frac{1}{2}(c+1)$ and $\beta = -\frac{1}{2}(s+1)$ will play an important role. Indeed, we will use the following as reference measures to generate the exponential families of distributions which are the central object of our study in this paper. These reference measures are

$$\mu_G(dx) = H_G(-\tfrac{1}{2}(c+1), -\tfrac{1}{2}(s+1); x)\mathbf{1}_{Q_G}(x)\,dx, \tag{2.6}$$

$$\nu_G(dy) = H_G(\tfrac{1}{2}(c+1), \tfrac{1}{2}(s+1); \varphi(y))\mathbf{1}_{P_G}(y)\,dy. \tag{2.7}$$

Applying (2.4), we see that $\nu_G$ is the image of $\mu_G$ by the mapping $x \mapsto y = \hat{x}^{-1}$ and that conversely $\mu_G$ is the image of $\nu_G$ by $y \mapsto x = \varphi(y)$.

Let $M_d^+$ denote the cone of positive definite matrices of order $d$ and $L(\mathbb{R}^p, \mathbb{R}^q)$ denote the space of linear transformations from $\mathbb{R}^p$ to $\mathbb{R}^q$. For $x \in Q_G$, $x_{C_j}, j = 1, \ldots, k$, are well defined and it will be convenient to use the following standard notation for various block submatrices:

$$
\begin{aligned}
x_{S_j} &= x_{<j>}, & x_{R_j, S_j} &= x_{[j>} = x_{<j]}^t, \\
x_{[j]} &= x_{R_j}, & x_{[j]\cdot} &= x_{[j]} - x_{[j>} x_{<j>}^{-1} x_{<j]},
\end{aligned}
\tag{2.8}
$$

where $x_{<j>} \in M_{s_j}^+, x_{[j]\cdot} \in M_{c_j - s_j}^+, x_{[j>} \in L(\mathbb{R}^{c_i - s_j}, \mathbb{R}^{s_j})$. It is understood here that $x_{[1]} = x_{[1]\cdot} = x_{C_1}$ whereas both $x_{<1>}$ and $x_{[1>}$ vanish. With this notation, we have, for example, $|\hat{x}| = \prod_{j=1}^k |x_{[j]\cdot}|$. In the proof of our main theorems, we will need to split the trace $\langle x, y \rangle$ for $x \in Q_G$ and $y \in P_G$ following a perfect order of the cliques as given in the following lemma.

LEMMA 2.2. *Let $G$ be a decomposable graph and let $C_1, \ldots, C_k$ be a perfect order of its cliques. For $x \in Q_G$ and $y \in P_G$ with $y = \hat{\sigma}^{-1}$ and $\sigma \in Q_G$, we have*

$$
\begin{aligned}
\langle x, y \rangle = \langle x, \hat{\sigma}^{-1} \rangle = \sum_{i=1}^k [&\langle x_{[i]\cdot}, \sigma_{[i]\cdot}^{-1} \rangle + \langle (x_{[i>} x_{<i>}^{-1} - \sigma_{[i>} \sigma_{<i>}^{-1}), \\
&\sigma_{[i]\cdot}^{-1} (x_{[i>} x_{<i>}^{-1} - \sigma_{[i>} \sigma_{<i>}^{-1}) x_{<i>} \rangle ].
\end{aligned}
\tag{2.9}
$$

This is a direct consequence of (2.3) and the following standard splitting of the trace for two positive definite matrices $u = \begin{pmatrix} u_1 & u_{12} \\ u_{21} & u_2 \end{pmatrix}$ and $v = \begin{pmatrix} v_1 & v_{12} \\ v_{21} & v_2 \end{pmatrix}$:

$$
\begin{aligned}
\langle u, v \rangle &= \langle u_1, v_{1\cdot 2} \rangle + \langle u_{2\cdot 1}, v_2 \rangle \\
&\quad + \langle (u_{21} u_1^{-1} + v_2^{-1} v_{21}), v_2 (u_{21} u_1^{-1} + v_2^{-1} v_{21}) u_1 \rangle,
\end{aligned}
\tag{2.10}
$$

and its corresponding expression if we write $v = \hat{\sigma}^{-1}$ with $\hat{\sigma}$ also positive definite,

$$
\begin{aligned}
\langle u, \hat{\sigma}^{-1} \rangle &= \langle u_1, \sigma_1^{-1} \rangle + \langle u_{2\cdot 1}, \sigma_{2\cdot 1}^{-1} \rangle \\
&\quad + \langle (u_{21} u_1^{-1} - \sigma_{21} \sigma_1^{-1}), \sigma_{2\cdot 1}^{-1} (u_{21} u_1^{-1} - \sigma_{21} \sigma_1^{-1}) u_1 \rangle.
\end{aligned}
\tag{2.11}
$$

We also recall the following basic results that will be used throughout our proofs.

LEMMA 2.3 ([13]). *The Jacobian of the change of variables*

$$
x \in Q_G \mapsto y = (x_{C_1}, x_{[i]\cdot}, x_{[i>} x_{<i>}^{-1}, i = 2, \ldots, k)
\tag{2.12}
$$

*is*

$$
\left| \frac{dy}{dx} \right| = \prod_{j=2}^k |x_{<j>}|^{c_j - s_j}.
\tag{2.13}
$$

The following lemma gives a Gaussian distribution we shall use later.

LEMMA 2.4. *For $x$ and $\sigma$ in $Q_G$, and for $L = L(R^{c_i-s_i}, R^{s_i})$ we have*

$$(2.14) \quad \int_L e^{-\langle (x_{[i>}x^{-1}_{<i>} - \sigma_{[i>}\sigma^{-1}_{<i>}), \sigma^{-1}_{[i].}(x_{[i>}x^{-1}_{<i>} - \sigma_{[i>}\sigma^{-1}_{<i>})x_{<i>}\rangle} \, d(x_{[i>}x^{-1}_{<i>})$$

$$= \pi^{(c_i-s_i)s_i/2} \frac{|\sigma_{[i].}|^{s_i/2}}{|x_{<i>}|^{(c_i-s_i)/2}}.$$

The proof follows immediately from Theorem 3.1.1 in [14] by replacing $C, D, Y$ and $M$ in that theorem by $\sigma_{[i].}, x^{-1}_{<i>}, x_{[i>}x^{-1}_{<i>}$ and $\sigma_{[i>}\sigma^{-1}_{<i>}$, respectively.

Let us finally recall the definition of the multivariate Gamma function. For $p > \frac{r-1}{2}$, the $r$-multivariate Gamma function is

$$(2.15) \quad \Gamma_r(p) = \pi^{(1/4)r(r-1)} \prod_{j=1}^{r} \Gamma(p - \tfrac{1}{2}(j-1)).$$

In the sequel we will need the following two formulas which link multivariate gamma functions of different dimensions. For $c$ and $s$ two positive integers with $s < c$ and for $\alpha > \frac{c-1}{2}$ a real number, we have

$$(2.16) \quad \pi^{(c-s)s/2} \Gamma_{c-s}\left(\alpha - \frac{s}{2}\right) = \frac{\Gamma_c(\alpha)}{\Gamma_s(\alpha)},$$

$$(2.17) \quad \pi^{(c-s)s/2} \Gamma_{c-s}(\alpha) = \frac{\Gamma_c(\alpha)}{\Gamma_s(\alpha - (c-s)/2)}.$$

2.2. *Tools for homogeneous graphs.* In this subsection, we study some properties of homogeneous graphs.

DEFINITION 2.1. A graph $G$ is said to be homogeneous if it is decomposable and does not contain the graph $\overset{1}{\bullet} - \overset{2}{\bullet} - \overset{3}{\bullet} - \overset{4}{\bullet}$, called $A_4$, as an induced subgraph.

We will see in Theorem 2.2 below why such a graph is called homogeneous. We now need to introduce a number of concepts about undirected graphs.

DEFINITION 2.2. Given an undirected graph $G = (V, E)$, the associated digraph is the directed graph $G' = (V, E')$ with $E'$ derived from $\mathcal{E}$ by the following process. If $i, j \in V$, then the directed edge $(i, j)$ is in $E'$ if and only if

$$(2.18) \quad \{i\} \cup nb(i) \supseteq \{j\} \cup nb(j),$$

where $nb(i) = \{j;\ j \neq i, i \sim j\}$.

Note that $E'$ contains all $(i,i)$ for $i \in V$. We write $i \to j$ if and only if $(i,j) \in E'$. An edge in $G$ can either disappear in $G'$ or become directed or become bi-directed. Note that if $i \sim j$ in $G$ then $i \not\to j$ if and only if there exists $k, k \neq i, k \neq j$, such that the subgraph of $G$ induced by $\{i,j,k\}$ is $\overset{k}{\bullet} - \overset{j}{\bullet} - \overset{i}{\bullet}$. Note also that if $k \not\sim i$, then it is impossible to have both $i \to j$ and $k \to j$. In other words, the configuration $\overset{i}{\bullet} \to \overset{j}{\bullet} \leftarrow \overset{k}{\bullet}$ in $G'$ is forbidden. Here are two simple examples of digraphs associated to given graphs. For the sake of clarity, the loops $i \to i$ are not drawn on the digraph $G'$.

EXAMPLES. The graph $\overset{1}{\bullet} - \overset{2}{\bullet}$ becomes the graph $G'$ $\overset{1}{\bullet} \leftrightarrow \overset{2}{\bullet}$. The graph $A_4$ becomes $\overset{1}{\bullet} \leftarrow \overset{2}{\bullet}$ $\overset{3}{\bullet} \to \overset{4}{\bullet}$.

It is easy to see from (2.18) that if $G$ is an undirected graph and $G'$ its associated digraph, then the relation $i \to j$ defined on $V$ is a preorder relation, that is $i \to j$ and $j \to k$ implies $i \to k$. Denote by $R$ the induced equivalence relation defined on $V$ by

$$iRj \Leftrightarrow i \to j \quad \text{and} \quad j \to i \Leftrightarrow \{i\} \cup nb(i) = \{j\} \cup nb(j).$$

Denote by $\bar{i}$ the equivalence class in $V/R$ containing $i \in V$ and denote by $\bar{i} \preceq \bar{j}$ the partial order relation on $V/R$ induced by the preorder $i \to j$. As usual, when dealing with partial order, the notation $\bar{i} \prec \bar{j}$ means $\bar{i} \preceq \bar{j}$ and $\bar{i} \neq \bar{j}$. We now introduce the Hasse diagram of $V/R$.

DEFINITION 2.3. The Hasse diagram of $G$ is the digraph with vertex set $V_H = V/R$ and with edge set $E_H$ such that an edge $(\bar{i}, \bar{j})$ is in $E_H$ if

$$\bar{i} \neq \bar{j}, \qquad \bar{i} \preceq \bar{j},$$

and

$$\bar{i} \preceq \bar{k} \preceq \bar{j} \qquad \text{implies either } \bar{k} = \bar{i} \text{ or } \bar{k} = \bar{j}.$$

If $(\bar{i}, \bar{j}) \in E_H$, we write $\bar{i} \to \bar{j}$. The knowledge of the Hasse diagram of $G$ is equivalent to the knowledge of the partial order relation on $V/R$. If $\bar{i} \to \bar{j}$ then $\bar{j}$ is a child of $\bar{i}$ and $\bar{i}$ is a parent of $\bar{j}$. If $i$ and $j$ are in $V$ it will be convenient to write $i \to j$ when the corresponding equivalence classes satisfy $\bar{i} \to \bar{j}$. Let us give an example of construction of a Hasse diagram.

In Figure 1, we give a graph $G$, its associated digraph $G'$ and the corresponding Hasse diagram. In $G'$, the loops $(i,i)$ are omitted. Since $3 \to 7$ and $7 \to 3$, $\{3,7\}$ is an equivalence class denoted by $\bar{3}$ while the 5 other vertices are alone in their equivalence class, which we denote by $i$ rather than $\bar{i}$ for simplicity. On this particular example, the Hasse diagram is a rooted tree associated to a partial order such that the root 1 is the minimum. The four cliques correspond to the endpoints of the tree: $C_{\bar{3}} = \{1,2,3,7\}$, $C_4 = \{1,2,4\}$, $C_5 = \{1,2,5\}$, $C_6 = \{1,6\}$. The two separators $S_1 = \{1,2\}$ with multiplicity 2 and $S_2 = \{1\}$ with multiplicity 1 correspond to the other vertices of the diagram. Note that the graph $G$ is homogeneous since it does not contain any $A_4$ as an induced subgraph. The fact that in this example the Hasse diagram is a rooted tree and the graph is homogeneous is not a coincidence since we have the following characterization theorem.

THEOREM 2.2. *Let $G = (V, E)$ be a connected graph and let $G' = (V, E')$ be its associated digraph. The following properties are equivalent:*

1. *$G$ is homogeneous.*
2. *If $i \sim j$ then either $i \to j$ or $j \to i$ in $G'$.*
3. *The Hasse diagram of the partially ordered set $(V/R, \preceq)$ is a rooted tree such that its root $\bar{1}$ is the minimal point of $V/R$ and such that the number of children of a vertex is never equal to one.*
4. *$P_G$ is a homogeneous cone (i.e., its automorphism group acts on it transitively).*
5. *$Q_G$ is a homogeneous cone.*

We shall only use equivalences between 1, 2 and 3, which are easy to prove. The equivalence with 4 and 5 is stated for the curiosity of the reader. The homogeneous graphs are specially simple to handle. We call $T$ the set of vertices of the corresponding Hasse tree, so $T = V/R$. Consider the subset of $V$

$$V_{\bar{i}} = \bigcup_{\bar{j} \preceq \bar{i}} \bar{j}.$$

We gather the properties of the Hasse tree of $G$ in the following proposition.

PROPOSITION 2.2. *If $T$ is the Hasse rooted tree of a homogeneous graph $G$ with $k$ cliques and $k'$ minimal separators, we have that:*

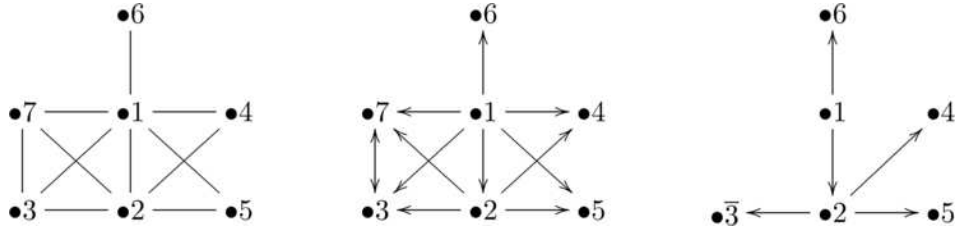

FIG. 1. *$G$, $G'$ and the Hasse diagram.*

1. *The mapping $\bar{i} \mapsto V_{\bar{i}}$, where $\bar{i} \in T$, gives a one to one correspondence between the cliques and minimal separators of $G$ and, respectively, the endpoints and non-endpoints of $T$. In particular, if $k > 1$ the root $\bar{1}$ is a minimal separator which is contained in all minimal separators and cliques of $G$ and the total number of vertices in $T$ is equal to $k + k'$.*

2. *All orders of the cliques are perfect. The multiplicity $\nu(V_s)$ of a separator $V_s$ is equal to the number of children of $s$ minus one.*

PROOF. 1. If $\bar{i} \in T$ then we observe that $V_{\bar{i}}$ is complete since if $j$ and $l$ are in $V_{\bar{i}}$ then either $\bar{j} \preceq \bar{l} \preceq \bar{i}$ or $\bar{l} \preceq \bar{j} \preceq \bar{i}$. In both cases $j \sim l$. Conversely, if $C \subset V$ is complete then $\overline{C} = \bigcup \{\bar{j} \in T; j \in C\}$ is contained in some $V_{\bar{i}}$. If not there exist $\bar{j}$ and $\bar{l}$ in $\overline{C}$ which are not comparable in the poset $T$ and therefore $j \not\sim l$, which contradicts the fact that $C$ is complete. Thus the maximal cliques are the $V_{\bar{i}}$'s where $\bar{i} \in T$ has no children, that is, $\bar{i}$ is an endpoint. Finally, if $\bar{i} \in T$ has children $\bar{j}$ and $\bar{l}$ then $V_{\bar{i}}$ is a minimal separator of $j$ and $l$ as can easily be seen. Conversely, if $j$ and $l$ are in $V$ with $j \not\sim l$ there exists a unique minimal separator between them which is $V_{\bar{i}}$ where $\bar{i} = \max\{s \in T;\ s \preceq \bar{j},\ s \preceq \bar{l}\}$.

2. Consider any order $(t(1), \ldots, t(k))$ of the endpoints of the tree and the corresponding order $(V_{t(1)}, \ldots, V_{t(k)})$ of the cliques. For $j = 2, \ldots, k$ and for $l = 1, \ldots, j-1$

$$s(l) = \max\{s \in T; s \prec t(j), s \prec t(l)\}.$$

Since $\bar{1} \preceq s(l) \prec t(j)$ for all $l = 1, \ldots, j-1$, then $s(l_j) = \max\{s(l); l = 1, \ldots, j-1\}$ exists and

$$V_{s(l_j)} = (V_{t(1)} \cup \cdots \cup V_{t(j-1)}) \cap V_{t(j)}$$

is a minimal separator contained in the clique $V_{t(l_j)}$ with $l_j < j$. Thus the order is perfect.

Now, given a minimal separator $V_s$, we show that the number $\nu(V_s)$ of $j$ such that there exists $l_j$ with $1 \leq l_j < j \leq k$ and $s = s(l_j)$, where $(l_j, s(l_j))$ is as defined above, is equal to $c(s) - 1$ where $c = c(s)$ is the number of children of $s$. Suppose first that $\nu(V_s) \geq c$. Then there exist endpoints $t(j_1), \ldots, t(j_c)$ of $T$ such that $j_1 < \cdots < j_c$ and such that $s = s(l_{j_1}) = \cdots = s(l_{j_c})$. Thus $s \prec t(j_1), \ldots, s \prec t(j_c)$. Furthermore $l_{j_1} < j_1$ and $s \prec t(l_{j_1})$. This implies that $s$ has at least $c+1$ children, a contradiction. Thus $\nu(V_s) \leq c(s) - 1$. Finally, one sees by induction that the number of edges of an undirected tree is the number of vertices minus one. Since $\sum_s c(s)$ is equal to the number of edges in the graph, this implies that $\sum_s (c(s) - 1) = k - 1$, where the sum is taken over the non-endpoints $s$ of $T$. To conclude the proof, we use the fact that by definition of the multiplicity of a minimal separator, the sum of the $\nu(V_s)$ is also $k - 1$. Thus $\sum_s [c(s) - 1 - \nu(V_s)] = (k-1) - (k-1) = 0$.

Since we have a null sum of nonnegative terms we get $\nu(V_s) = c(s) - 1$ for all minimal separators. $\square$

It follows from the proposition above that there is a one to one correspondence between the set of homogeneous graphs and the set of rooted trees with vertices weighted by positive integers and such that no vertex has exactly one child. Note that a complete graph is homogeneous. It is characterized by the fact that its Hasse diagram is just a point. A decomposable graph with only one separator is homogeneous. Its Hasse tree looks like a daisy. An undirected tree is decomposable but is not homogeneous in general. Finally it is possible to prove that if all orders of the cliques of a decomposable graph $G$ are perfect then $G$ is homogeneous.

**3. The Wishart families of Types I and II.** In this section, we define two families of Wishart distributions. We will study special cases in Section 3.2, the homogeneous case in Section 3.3 and the nonhomogeneous case in Section 3.4.

3.1. *Definitions.* Consider the two integrals

$$(3.1) \quad I(\alpha, \beta; y) = \int_{Q_G} e^{-\langle x, y \rangle} H_G(\alpha, \beta; x) \mu_G(dx) \qquad \text{for } y \in P_G,$$

$$(3.2) \quad J(\alpha, \beta; x) = \int_{P_G} e^{-\langle x, y \rangle} H_G(\alpha, \beta; \varphi(y)) \nu_G(dy) \qquad \text{for } x \in Q_G.$$

We define $\mathcal{A}$ to be the set of $(\alpha, \beta)$ such that $I(\alpha, \beta; y)$ converges for all $y \in P_G$ and such that $y \mapsto \frac{I(\alpha, \beta; y)}{H_G(\alpha, \beta; \varphi(y))}$ is a constant on $P_G$. This constant is a function on $\mathcal{A}$ that we denote by $\Gamma_{\mathrm{I}}(\alpha, \beta)$. Similarly we define $\mathcal{B}$ to be the set of $(\alpha, \beta)$ such that $J(\alpha, \beta; x)$ converges for all $x \in Q_G$ and such that $x \mapsto \frac{J(\alpha, \beta; x)}{H_G(\alpha, \beta; x)}$ is a constant on $Q_G$. This constant is a function on $\mathcal{B}$ that we denote by $\Gamma_{\mathrm{II}}(\alpha, \beta)$. The sets $\mathcal{A}$ and $\mathcal{B}$ will be studied in Sections 3.3, 3.4 and 3.5.

We note here that since $\mu_G(du)$ is the image of the measure $\nu_G(dy)$ under the mapping $y \mapsto u = \varphi(y)$ (see Section 2.1), (3.2) can be written

$$J(\alpha, \beta; x) = \int_{Q_G} e^{-\langle x, \hat{u}^{-1} \rangle} H_G(\alpha, \beta; u) \mu_G(du) \qquad \text{for } x \in Q_G.$$

This expression of (3.2) and the passage from $y \in P_G$ to $u = \varphi(y) \in Q_G$ will be used several times in the remainder of the paper for defining the inverse Type II Wishart and to perform various computations. The Wishart distributions of Type I will be the probabilities

$$\frac{1}{\Gamma_{\mathrm{I}}(\alpha, \beta) H_G(\alpha, \beta; \varphi(y))} e^{-\langle x, y \rangle} H_G(\alpha, \beta; x) \mu_G(dx),$$

defined on $Q_G$ and indexed by the parameters $(\alpha, \beta; y)$ in $\mathcal{A} \times P_G$. To follow the standard notation for distributions related to the Wishart, when $y \in P_G$ is the parameter of the Type I Wishart, we often write $y = \hat{\sigma}^{-1}$ with $\sigma \in Q_G$ so that, for $\sigma \in Q_G, (\alpha, \beta) \in \mathcal{A}$, the Type I Wishart distribution can be written

$$(3.3) \qquad W_{Q_G}(\alpha, \beta, \sigma; dx) = e^{-\langle x, \hat{\sigma}^{-1}\rangle} \frac{H_G(\alpha, \beta; x)}{\Gamma_I(\alpha, \beta) H_G(\alpha, \beta; \sigma)} \mu_G(dx).$$

The Wishart distributions of Type II will be the probabilities

$$(3.4) \qquad W_{P_G}(\alpha, \beta, \theta; dy) = e^{-\langle \theta, y\rangle} \frac{H_G(\alpha, \beta; \varphi(y))}{\Gamma_{II}(\alpha, \beta) H_G(\alpha, \beta; \theta)} \nu_G(dy)$$

defined on $P_G$ and indexed by the parameters $(\alpha, \beta; \theta)$ in $\mathcal{B} \times Q_G$. We therefore consider the following two natural exponential families.

DEFINITION 3.1. For $(\alpha, \beta) \in \mathcal{A}$, the Type I Wishart family of distributions is defined by

$$(3.5) \qquad \mathcal{F}_{(\alpha,\beta),I} = \{W_{Q_G}(\alpha, \beta, \sigma; dx), \sigma \in Q_G\}.$$

DEFINITION 3.2. For $(\alpha, \beta) \in \mathcal{B}$, the Type II Wishart family of distributions is defined by

$$(3.6) \qquad \mathcal{F}_{(\alpha,\beta),II} = \{W_{P_G}(\alpha, \beta, \theta; dy), \theta \in Q_G\}.$$

Following the pattern of what is done for the Wishart distribution, we now define Type I and Type II inverse Wishart and $F$ distributions.

DEFINITION 3.3. Let $G$ be given. If $X \sim W_{Q_G}(\alpha, \beta, \sigma)$ where $(\alpha, \beta) \in \mathcal{A}$ and $\sigma \in Q_G$, then $Y = \hat{X}^{-1}$ is said to follow the inverse Type I Wishart, defined on $P_G$, and its distribution is

$$(3.7) \qquad IW_{Q_G}(\alpha, \beta, \sigma; dy) = \frac{e^{-\langle \varphi(y), \hat{\sigma}^{-1}\rangle} H_G(\alpha, \beta; \varphi(y))}{\Gamma_I(\alpha, \beta) H_G(\alpha, \beta; \sigma)} \nu_G(dy).$$

The distribution (3.7) is clearly immediately derived from the distribution (3.3) by recalling that $x = \varphi(y)$ and that $\nu_G(dy)$ is the image of $\mu_G(dx)$ by the mapping $x \mapsto y = \hat{x}^{-1}$.

DEFINITION 3.4. Let $G$ be given. If $Y \sim W_{P_G}(\alpha, \beta, \theta)$ where $(\alpha, \beta) \in \mathcal{B}$ and $\theta \in Q_G$, then $X = \varphi(Y)$ is said to follow the inverse Type II Wishart, defined on $Q_G$, and its distribution is

$$(3.8) \qquad IW_{P_G}(\alpha, \beta, \theta; dx) = \frac{e^{-\langle \theta, \hat{x}^{-1}\rangle} H_G(\alpha, \beta; x)}{\Gamma_{II}(\alpha, \beta) H_G(\alpha, \beta; \theta)} \mu_G(dx).$$

Here too, the density (3.8) is immediately derived from (3.4).

Let $\mathcal{B} - \mathcal{A} = \{(\alpha' - \alpha, \beta' - \beta) : (\alpha', \beta') \in \mathcal{B}, (\alpha, \beta) \in \mathcal{A}\}$. Since $\mathcal{B} - \mathcal{A} \subset \mathcal{B}$ and $\mathcal{A} - \mathcal{B} \subset \mathcal{A}$ are false in general, as will be seen, for example, when $G = A_4$, to give the following definition of the $F$ distributions, we will have to insure that the parameters $\alpha' - \alpha$ and $\beta' - \beta$ are in the correct sets.

DEFINITION 3.5. Let $\theta$ and $\sigma$ be in $P_G$ and $Q_G$, and let $(\alpha, \beta) \in \mathcal{A}$ and $(\alpha', \beta') \in \mathcal{B}$. Then

1. for $(\alpha' - \alpha, \beta' - \beta) \in \mathcal{B}, \sigma \in Q_G$, the $F$ distribution of the first kind with parameters $(\alpha, \beta, \alpha', \beta', \sigma)$ is the distribution on $Q_G$

$$\frac{\Gamma_{\mathrm{II}}(\alpha' - \alpha, \beta' - \beta)}{\Gamma_{\mathrm{I}}(\alpha, \beta)\Gamma_{\mathrm{II}}(\alpha', \beta')} H_G(-\alpha', -\beta'; \sigma)$$
$$\times H_G(\alpha' - \alpha, \beta' - \beta; \sigma + x) H_G(\alpha, \beta; x) \mu_G(dx);$$

2. for $(\alpha - \alpha', \beta - \beta') \in \mathcal{A}, \theta \in Q_G$, the $F$ distribution of the second kind with parameters $(\alpha, \beta, \alpha', \beta', \theta)$ is the distribution on $P_G$

$$\frac{\Gamma_{\mathrm{I}}(\alpha - \alpha', \beta - \beta')}{\Gamma_{\mathrm{I}}(\alpha, \beta)\Gamma_{\mathrm{II}}(\alpha', \beta')} H_G(-\alpha, -\beta; \varphi(\theta))$$
$$\times H_G(\alpha - \alpha', \beta - \beta'; \varphi(\theta + y)) H_G(\alpha', \beta'; \varphi(y)) \nu_G(dy).$$

Note here again that the lack of multiplicative structure on $P_G$ and $Q_G$ when $G$ is not complete prevents us from relating these distributions to some form of quotient $X/X'$ of independent random variables with distributions $W_{Q_G}(\alpha, \beta; \sigma)$ and $W_{Q_G}(\alpha', \beta'; \sigma)$, respectively. A study of the multivariate $F$ distribution when $G$ is complete can be found in [15]. We could also define rather explicitly Beta distributions of Type I by introducing the conditional distributions $X|X + X'$, where $X \sim W_{Q_G}(\alpha, \beta, \sigma)$ and $X' \sim W_{Q_G}(\alpha', \beta', \sigma)$ are independent such that $(\alpha + \alpha', \beta + \beta')$ is still in $\mathcal{A}$. Again, since the cone $Q_G$ has no special multiplicative structure, these Beta distributions unfortunately do not seem to enjoy properties linking them to some ratio analogous to $X/(X + X')$ as happens when the graph is complete. The same problem arises with Beta distributions of Type II. Finally, we could also consider the distribution

$$\frac{e^{-\langle x, y \rangle} H_G(\alpha, \beta; x)}{I(\alpha, \beta; y)} \mu_G(dx)$$

where we only require that $(\alpha, \beta)$ be such that $I(\alpha, \beta; y)$ defined by (3.1) converges. Under such generality, these distributions have no interesting properties: their Laplace transforms are not explicit, their family is not stable by convolution as our Wishart distributions of Type I are (see Proposition 3.2) and they have no hyper Markov property. A similar remark holds for $J$ defined by (3.2) and Type II Wisharts.

3.2. *The hyper and inverse hyper inverse Wishart distributions.* We first observe that when $G$ is complete, both Type I and Type II Wishart distributions coincide with the ordinary Wishart distribution. We will see now that for special values of $(\alpha, \beta)$, the Type I and II Wisharts are, respectively, the hyper Wishart as defined by Dawid and Lauritzen [7] and the $G$-Wishart first identified by Roverato [16] as the inverse of the hyper inverse Wishart defined also by Dawid and Lauritzen [7]. To describe these distributions, it is convenient to fix a perfect order of the cliques.

*The hyper Wishart on $Q_G$.* Let $G$ be given and let $p$ be a scalar. Let $\mathcal{A}_1$ be the one-dimensional subset of $\mathbb{R}^{k+k'}$ defined as

$$\mathcal{A}_1 = \left\{ (\alpha, \beta) | \alpha(C) = p, C \in \mathcal{C}, \beta(S) = p, S \in \mathcal{S} \text{ with } p > \max_{C \in \mathcal{C}} \tfrac{1}{2}(|C| - 1) \right\}.$$

For $(\alpha, \beta) \in \mathcal{A}_1$ we then have

$$(3.9) \qquad W_{Q_G}(\alpha, \beta, \sigma; dx) \propto \frac{\prod_{i=1}^{k} w_{c_i}(p, \sigma_{C_i}; x_{C_i})}{\prod_{i=2}^{k} w_{s_i}(p, \sigma_{S_i}; x_{S_i})} \mathbf{1}_{Q_G}(x) \, dx$$

with

$$w_{c_i}(p, \sigma_{C_i}; x_{C_i}) = \frac{|x_{C_i}|^{p-(c_i+1)/2}}{\Gamma_{c_i}(p) |\sigma_{C_i}|^p} e^{-\langle x_{C_i}, \sigma_{C_i}^{-1} \rangle}.$$

We note that the expression of $W_{Q_G}(\alpha, \beta, \sigma; dx)$ in (3.9) does not depend on the chosen perfect order of the cliques. The expression on the right-hand side of (3.9) is a Markov combination of Wishart distributions with shape parameters $p$ and scale parameters $\sigma_{C_i}$ and $\sigma_{S_i}$ on the cliques and separators of $G$, respectively. By Theorem 2.6 of [7], it is a distribution. It is in fact the hyper Wishart distribution, as defined in that same paper. Therefore both sides of (3.9) are equal and equal to the density of the hyper Wishart distribution and it follows immediately that $\mathcal{A}_1 \subset \mathcal{A}$ for any given $G$ and

$$\Gamma_{\mathrm{I}}(\alpha, \beta) = \frac{\prod_{i=1}^{k} \Gamma_{c_i}(p)}{\prod_{i=2}^{k} \Gamma_{s_i}(p)}.$$

*The $G$-Wishart on $P_G$.* Let $G$ be given and let $\delta > 0$ be a scalar. Let $\mathcal{B}_1$ be the one-dimensional subset of $\mathbb{R}^{k+k'}$ defined as

$$\mathcal{B}_1 = \{ (\alpha, \beta) | \alpha(C) = -\tfrac{1}{2}(\delta + |C| - 1), C \in \mathcal{C},$$
$$\beta(S) = -\tfrac{1}{2}(\delta + |S| - 1), S \in \mathcal{S}, \delta > 0 \}.$$

For $(\alpha, \beta) \in \mathcal{B}_1$ and for $x = \varphi(y)$

$$H_G(\alpha, \beta; \varphi(y)) \nu_G(dy) = \frac{\prod_{i=1}^{k} |x_{C_i}|^{-(\delta+c_i-1)/2+(c_i+1)/2}}{\prod_{i=2}^{k} |x_{S_i}|^{-(\delta+s_i-1)/2+(s_i+1)/2}} \mathbf{1}_{P_G}(y) \, dy$$

$$= \frac{\prod_{i=1}^{k} |x_{C_i}|^{-(\delta-2)/2}}{\prod_{i=2}^{k} |x_{S_i}|^{-(\delta-2)/2}} \mathbf{1}_{P_G}(y) \, dy$$

$$= |y|^{(\delta-2)/2} \mathbf{1}_{P_G}(y) \, dy,$$

where, as before, the expression of $H_G(\alpha, \beta; \varphi(y))\nu_G(dy)$ does not depend on any chosen perfect order of the cliques. Therefore, $W_{P_G}(\alpha, \beta, \theta; dy) \propto |y|^{(\delta-2)/2} e^{-\langle \theta, y \rangle} \, dy$ is the $G$-Wishart distribution first identified by Roverato [16] as the inverse of the hyper inverse Wishart. It follows immediately that $\mathcal{B}_1 \subset \mathcal{B}$, that

$$\Gamma_{\mathrm{II}}(\alpha, \beta) = \frac{\prod_{i=1}^{k} \Gamma_{c_i}((\delta + c_i - 1)/2)}{\prod_{i=2}^{k} \Gamma_{s_i}((\delta + s_i - 1)/2)}$$

and that the Type II Wishart is the $G$-Wishart defined on $P_G$ for $\delta > 0$, $\theta \in Q_G$. The distribution of this special Type II Wishart is

$$W_{P_G}(\alpha, \beta, \theta) = \frac{\prod_{i=1}^{k} |\theta_{C_i}|^{(\delta+c_i-1)/2}}{\prod_{i=1}^{k} \Gamma_{c_i}((\delta + c_i - 1)/2)} \frac{\prod_{i=2}^{k} \Gamma_{s_i}((\delta + s_i - 1)/2)}{\prod_{i=2}^{k} |\theta_{S_i}|^{(\delta+s_i-1)/2}}$$
$$\times |y|^{(\delta-2)/2} e^{-\langle \theta, y \rangle} \mathbf{1}_{P_G}(y) \, dy.$$

3.3. *The homogeneous case.* We now consider the Type I and II Wishart distributions when the graph $G$ is homogeneous as defined in Section 2.2. Our aim is to identify the sets $\mathcal{A}$ and $\mathcal{B}$ and the values of the normalizing constants $\Gamma_{\mathrm{I}}(\alpha, \beta)$ and $\Gamma_{\mathrm{II}}(\alpha, \beta)$. It is convenient to introduce the following notation, consistent with the notation introduced in (2.8). For $u \in V/R$, we define

$$x_{[u]} = (x_{ij})_{\bar{i}=\bar{j}=u}, \qquad x_{[u>} = (x_{ij})_{\bar{i}=u, \bar{j} \prec u},$$
$$x_{<u>} = (x_{ij})_{\bar{i} \prec u, \bar{j} \prec u}, \qquad x_{[u]\cdot} = x_{[u]} - x_{[u>} x_{<u>}^{-1} x_{<u]}.$$

We also adopt the convention that a vertex of the Hasse tree $(T, E_H)$ will be denoted by $t$ if it is an endpoint of the tree and by $q$ if it is not an endpoint. From Proposition 2.2, to each $t$ corresponds a unique clique $C_t = \bigcup_{u \in T, u \preceq t}[u]$ and therefore a number $\alpha_t = \alpha(C_t)$. And, to each $q$ corresponds a unique minimal separator $S_q = \bigcup_{u \in T, u \preceq q}[u]$ and therefore a number $\beta_q = \beta(S_q)$. The positive integer $\nu(q)$ is the number of children of $q$ minus one. From Proposition 2.2, this is also the multiplicity of $S_q$. With these conventions, for each $u \in T$, we write

$$\rho_u = \rho_u(\alpha, \beta) = \sum_{u \preceq t} \alpha_t - \sum_{u \preceq q} \nu(q) \beta_q.$$

We define $n_u$ to be the cardinality of the vertex $u$ of the Hasse tree, that is, the number of vertices in $V$ that are in the vertex $u$ of the Hasse tree of $G$.

We also define

$$m_u = n_u + \sum_{v \prec u} n_v = \sum_{v \preceq u} n_v. \tag{3.10}$$

The following two theorems give $\mathcal{A}, \mathcal{B}$ and the corresponding normalizing constants for the Wisharts of Type I and II in the homogeneous case.

THEOREM 3.1. *Let $G$ be a homogeneous graph. Then*

$$\mathcal{A} = \left\{ (\alpha, \beta) \Big| \rho_u > \tfrac{1}{2}\left(\sum_{v \preceq u} n_v - 1\right), u \in T \right\}.$$

*More specifically for $(\alpha, \beta) \in \mathcal{A}$ and $\sigma \in Q_G$, the integral* (3.1) *converges and*

$$\int_{Q_G} e^{-\langle x, \hat{\sigma}^{-1} \rangle} H_G(\alpha, \beta; x) \mu_G(dx)$$

$$= \prod_{u \in T} \pi^{\sum_{v \prec u} n_u n_v / 2} |\sigma_{[u]}|^{\rho_u} \Gamma_{n_u}\left(\rho_u - \sum_{v \prec u} \frac{n_v}{2}\right) \tag{3.11}$$

$$= H_G(\alpha, \beta; \sigma) \prod_{u \in T} \pi^{\sum_{v \prec u} n_u n_v / 2} \Gamma_{n_u}\left(\rho_u - \sum_{v \prec u} \frac{n_v}{2}\right).$$

THEOREM 3.2. *Let $G$ be a homogeneous graph. Then*

$$\mathcal{B} = \left\{ (\alpha, \beta) \Big| -\rho_u > \tfrac{1}{2}\left(\sum_{v \succeq u} n_v - 1\right), u \in T \right\}.$$

*More specifically for $(\alpha, \beta) \in \mathcal{B}$ and $\theta \in Q_G$, the integral* (3.2) *converges and*

$$\int_{P_G} e^{-\langle y, \theta \rangle} H_G(\alpha, \beta; \varphi(y)) \nu_G(dy)$$

$$= \prod_{u \in T} \pi^{\sum_{v \prec u} n_u n_v / 2} |\theta_{[u]}|^{-\rho_u} \Gamma_{n_u}\left(-\rho_u - \sum_{v \succ u} n_v / 2\right) \tag{3.12}$$

$$= H_G(\alpha, \beta; \theta) \prod_{u \in T} \pi^{\sum_{v \prec u} n_u n_v / 2} \Gamma_{n_u}\left(-\rho_u - \sum_{v \succ u} \frac{n_v}{2}\right).$$

We note that the parameter sets $\mathcal{A}$ and $\mathcal{B}$ in the homogeneous case are $(k + k')$-dimensional. The proof of Theorem 3.1 follows the same line as that of the proof of Theorem 3.3 given in the Appendix. It is based on Proposition 3.1 given below and on the following analog of formula (2.9) for the traces:

for $x$ and $\sigma$ in $Q_G$,

$$\langle x, \hat{\sigma}^{-1}\rangle = \sum_{u\in T}[\langle x_{[u]\cdot}, \sigma_{[u]\cdot}^{-1}\rangle$$
(3.13)
$$+ \langle (x_{[u>}x_{<u>}^{-1} - \sigma_{[u\rangle}\sigma_{<u>}^{-1}),$$
$$\sigma_{[u]\cdot}^{-1}(x_{[u>}x_{<u>}^{-1} - \sigma_{[u>}\sigma_{<u>}^{-1})x_{<u>}\rangle],$$

where it is understood that, as in (2.8) and (2.9), for $u = \overline{1}$, the root of $T$, the summand reduces to $\langle x_1, \sigma_1^{-1}\rangle$. Then, using Proposition 3.1 and formula (3.13), the integral in (3.11) is obtained by a series of standard integrations. The proof of (3.13) is parallel to the proof of (2.9) and will not be given here. The proof of Proposition 3.1 is given in the Appendix.

PROPOSITION 3.1. *For $G$ homogeneous, the image of $H_G(\alpha, \beta; x)\mu_G(dx)$ under the mapping*

(3.14) $$x = (x_{[u]}, x_{[u>}, u \in T) \mapsto (x_{[u]\cdot}, x_{[u>}x_{<u>}^{-1}, u \in T)$$

*is*

$$H_G^*(\alpha, \beta; x)\mu_G^*(dx_{[u]\cdot}, dx_{[u>}x_{<u>}^{-1}, u \in T)$$
$$= \prod_{u\in T}|x_{[u]\cdot}|^{\lambda_u-(n_u+1)/2}\,dx_{[u]\cdot}\,d(x_{[u>}x_{<u>}^{-1}),$$

*where*

$$\lambda_u = \rho_u + \sum_{v\succ u}\frac{n_v}{2} - \sum_{v\prec u}\frac{n_v}{2}.$$

The proof of Theorem 3.2 also follows the general lines of the proof of Theorem 3.4 given in the Appendix. We first observe that the image of $H_G(\alpha, \beta; \varphi(y)) \times \nu_G(dy)$ under the change of variable $y \mapsto x = \varphi(y)$ is $H_G(\alpha, \beta; x)\mu_G(dx)$ so that

$$\int_{P_G} e^{-\langle y, \theta\rangle}H_G(\alpha, \beta; \varphi(y))\nu_G(dy)$$
$$= \int_{Q_G} e^{-\langle \theta, \hat{x}^{-1}\rangle}\frac{\prod_{j=1}^{k}|x_{C_j}|^{\alpha_j-(c_j+1)/2}}{\prod_{j=2}^{k}|x_{S_j}|^{\beta_j-(s_j+1)/2}}\,dx,$$

where, as usual, the integral on the right-hand side of the equation above does not depend upon the chosen perfect order of the cliques. We then use (3.14) and (3.13) applied to $\langle \theta, \hat{x}^{-1}\rangle$ to obtain the expression of the integral in (3.12) by a series of standard integrations.

Using Proposition 3.1, (3.11) and (3.13) it is fairly straightforward to show that the image of the Type I Wishart by the change of variable (3.14) is the distribution

$$W^*_{Q_G}(\alpha, \beta, \sigma; dx_{[u]\cdot}, d(x_{[u>}x^{-1}_{<u>})), u \in T)$$

$$= \prod_{u \in T} \left[ |x_{[u]\cdot}|^{\lambda_u - (n_u+1)/2} e^{-\langle x_{[u]\cdot}, \sigma^{-1}_{[u]\cdot} \rangle} \right.$$

(3.15)
$$\times e^{-\langle (x_{[u>}x^{-1}_{<u>} - \sigma_{[u>}\sigma^{-1}_{<u>}), \sigma^{-1}_{[u]\cdot}(x_{[u>}x^{-1}_{<u>} - \sigma_{[u>}\sigma^{-1}_{<u>})x_{<u>} \rangle}$$

$$\times \left( \pi^{(1/2)n_u(\sum_{v \succ u} n_v)} |\sigma_{[u]\cdot}|^{\rho_u} \Gamma_{n_u}\left(\lambda_u - \tfrac{1}{2}\sum_{v \prec u} n_v\right) \right)^{-1}$$

$$\left. \times \mathbf{1}_{D_u}(x_{[u]\cdot}, (x_{[u>}x^{-1}_{<u>})) \, dx_{[u]\cdot} \, d(x_{[u>}x^{-1}_{<u>}) \right],$$

where $D_u = (M^+_{n_u} \times L(\mathbb{R}^{n_u}, \mathbb{R}^{m_u - n_u}))$. This distribution is exactly the distribution of the Wisharts defined by Andersson and Wojnar [3] on homogeneous cones. Therefore when $G$ is a homogeneous graph the Type I Wisharts coincide with the Wisharts of [3] for the homogeneous cone $Q_G$ corresponding to $G$. Since the dual $P_G$ of $Q_G$ is also homogeneous, we could also show that the Type II Wisharts correspond to the Wisharts as defined by Andersson and Wojnar [3] on $P_G$. However, there are many other homogeneous cones not of the form $P_G$ and $Q_G$. Our calculations are simpler and self contained in the particular cases that we investigate here.

Using (3.15), Proposition 3.1, (3.13) and (3.12), we obtain the image of the inverse Type II Wishart by the change of variable (3.14). The image of the distribution of $X \sim IW_{P_G}(\alpha, \beta, \theta)$, the inverse of the Type II Wishart when $G$ is homogeneous, is given by

$$IW^*_{P_G}(\alpha, \beta, \theta; dx_{[u]\cdot}, dx_{[u>}x^{-1}_{<u>}, u \in T)$$

$$= \prod_{u \in T} \left[ |x_{[u]\cdot}|^{\lambda_u - (n_u+1)/2} e^{-\langle x^{-1}_{[u]\cdot}, \theta^{-1}_{[u]\cdot} \rangle} \right.$$

(3.16)
$$\times e^{-\langle (x_{[u>}x^{-1}_{<u>} - \theta_{[u>}\theta^{-1}_{<u>}), \theta_{<u>}(x^{-1}_{<u>}x_{<u]} - \theta^{-1}_{<u>}\theta_{<u]})x^{-1}_{[u]\cdot} \rangle}$$

$$\times \left( \pi^{(1/2)n_u \sum_{v \prec u} n_v} |\theta_{[u]\cdot}|^{-\rho_u} \Gamma_{n_u}\left(-\rho_u - \sum_{v \succ u} \frac{n_v}{2}\right) \right)^{-1}$$

$$\left. \times \mathbf{1}_{D_u}(x_{[u]\cdot}, (x_{[u>}x^{-1}_{<u>})) \, dx_{[u]\cdot} \, d(x_{[u>}x^{-1}_{<u>}) \right].$$

EXAMPLE. Consider the graph $G_0$ (see Figure 2). We index each clique according to the vertex of the Hasse tree of $G_0$ which represents it. Thus $C_3 = \{1,2,3\}$, $C_4 = \{1,2,4\}$, $C_5 = \{1,2,5\}$, $C_6 = \{1,6\}$. Minimal separators are $S_1 = \{1\}$ and $S_2 = \{1,2\}$ with $\nu(S_2) = 2$ and $\nu(S_1) = 1$. We set

$$\alpha(C_3) = \alpha_3, \quad \alpha(C_4) = \alpha_4, \quad \alpha(C_5) = \alpha_5, \quad \alpha(C_6) = \alpha_6,$$
$$\beta(S_1) = \beta_1, \quad \beta(S_2) = \beta_2.$$

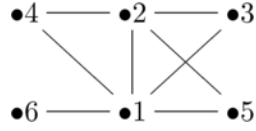

Fig. 2. *Graph $G_0$.*

The Hasse tree corresponding to $G_0$ is identical to the Hasse diagram of Figure 1 with $\overline{3}$ replaced by 3. Since the cardinality of all the vertices of the Hasse tree is 1, for the sake of simplicity we will denote the vertices of the tree by 1, 2, 3, 4, 5 and 6 so that we have

$$n_i = 1, \quad i = 1,\ldots,6 \quad \text{and} \quad 1 \preceq 6, \ 1 \preceq 2, \ 2 \preceq 3, \ 2 \preceq 4, \ 2 \preceq 5,$$
$$\rho_1 = \alpha_3 + \alpha_4 + \alpha_5 + \alpha_6 - \beta_1 - 2\beta_2, \quad \rho_2 = \alpha_3 + \alpha_4 + \alpha_5 - 2\beta_2,$$
$$\rho_3 = \alpha_3, \quad \rho_4 = \alpha_4, \quad \rho_5 = \alpha_5, \quad \rho_6 = \alpha_6,$$

and

$$\lambda_1 = \rho_1 + \tfrac{5}{2}, \quad \lambda_2 = \rho_2 + \tfrac{3}{2} - \tfrac{1}{2},$$
$$\lambda_3 = \rho_3 - \tfrac{2}{2}, \quad \lambda_4 = \rho_4 - \tfrac{2}{2},$$
$$\lambda_5 = \rho_5 - \tfrac{2}{2}, \quad \lambda_6 = \rho_6 - \tfrac{1}{2}.$$

Therefore,

$$\mathcal{A} = \{(\alpha,\beta) | \rho_1 > 0, \ \rho_2 > \tfrac{1}{2}, \ \rho_i > 1, i = 3,4,5, \ \rho_6 > \tfrac{1}{2}\}$$
(3.17)
$$= \{(\alpha,\beta) | \alpha_i > 1, \ i = 3,4,5, \ \alpha_6 > \tfrac{1}{2},$$
$$\alpha_3 + \alpha_4 + \alpha_5 + \alpha_6 - 2\beta_2 - \beta_1 > 0, \ \alpha_3 + \alpha_4 + \alpha_5 - 2\beta_2 > \tfrac{1}{2}\},$$
$$\mathcal{B} = \{(\alpha,\beta) | \rho_1 < -\tfrac{5}{2}, \rho_2 < -\tfrac{3}{2}, \rho_i < 0, i = 3,4,5,6\}$$
(3.18)
$$= \{(\alpha,\beta) | \alpha_i < 0, \ i = 3,4,5,6,$$
$$\alpha_3 + \alpha_4 + \alpha_5 + \alpha_6 - 2\beta_2 - \beta_1 < -\tfrac{5}{2}, \ \alpha_3 + \alpha_4 + \alpha_5 - 2\beta_2 < -\tfrac{3}{2}\}.$$

3.4. *The nonhomogeneous case.* We now consider a nonhomogeneous graph $G$, that is, a graph containing $A_4$ as an induced subgraph. As in the case of homogeneous graphs, our aim is to identify $\mathcal{A}, \mathcal{B}$ and the corresponding eigenvalues. We will see that we are, in fact, only able to identify a subset of $\mathcal{A}$ and $\mathcal{B}$ and the corresponding eigenvalues $\Gamma(\alpha, \beta)$. The results are given in Theorem 3.3 and Theorem 3.4 below. For $G$ a noncomplete decomposable graph, let $P = (C_1, \ldots, C_k)$ be a perfect order of the family $\mathcal{C}$ of its cliques and $(S_2, \ldots, S_k)$ be the associated sequence of minimal separators. Recall that $c_j = |C_j|$ and $s_j = |S_j|$ denote the cardinality of $C_j$ and $S_j$, respectively. For given $\alpha$ and $\beta$ we write $\alpha_j = \alpha(C_j)$ and $\beta_j = \beta(S_j)$. For a given minimal separator $S$ we write

$$J(P, S) = \{j = 2, \ldots, k | S_j = S\},$$

and for a given perfect order $P$ of the cliques, we define $A_P$ to be the set of $(\alpha, \beta)$ such that:

1. $\sum_{j \in J(P,S)} \alpha_j - \nu(S)\beta(S) = 0$, for all $S$ different of $S_2$;
2. $\alpha_j - \frac{c_j - 1}{2} > 0$ for all $C_j \in \mathcal{C}$;
3. $\alpha_1 + \delta_2 > \frac{s_2 - 1}{2}$ where $\delta_2 = \sum_{j \in J(P, S_2)} \alpha_j - \nu(S_2)\beta_2$.

Recall also that $\Gamma_n(p)$ is defined in (2.15). To avoid trivialities, in the following statements we assume that $G$ is not complete in order to have at least one minimal separator. Theorems 3.3 and 3.4 are useful only for nonhomogeneous graphs, since stronger results, Theorems 3.1 and 3.2, are available for homogeneous graphs.

THEOREM 3.3. *Let $G$ be a noncomplete decomposable graph and let $P$ be a perfect order of its cliques. Then $A_P \subset \mathcal{A}$. More specifically for $y \in P_G$ and for $(\alpha, \beta) \in A_P$ the integral (3.1) converges and*

$$(3.19) \quad \int_{Q_G} e^{-\langle x, y \rangle} H_G(\alpha, \beta; x) \mu_G(dx) = \Gamma_\mathrm{I}(\alpha, \beta) H_G(\alpha, \beta; \varphi(y)),$$

*where*

$$(3.20) \quad \Gamma_\mathrm{I}(\alpha, \beta) = \Gamma_{s_2}(\alpha_1 + \delta_2) \frac{\Gamma_{c_1}(\alpha_1)}{\Gamma_{s_2}(\alpha_1)} \prod_{q=2}^{k} \frac{\Gamma_{c_q}(\alpha_q)}{\Gamma_{s_q}(\alpha_q)}.$$

Equivalently, if we write $y = \hat{\sigma}^{-1}$ with $\sigma \in Q_G$, (3.19) can be rewritten as

$$(3.21) \quad \int_{Q_G} e^{-\langle x, \hat{\sigma}^{-1} \rangle} H_G(\alpha, \beta; x) \mu_G(dx) = \Gamma_\mathrm{I}(\alpha, \beta) H_G(\alpha, \beta; \sigma).$$

To study $\mathcal{B}$ for a nonhomogeneous graph and give the normalizing constant of the Type II Wishart, we now need to define, for a given $P$, the set $B_P$ to be the set of $(\alpha, \beta)$ such that

1. $\sum_{j \in J(P,S)} (\alpha_j + \frac{1}{2}(c_j - s_j)) - \nu(S)\beta(S) = 0,$ for all $S$ different from $S_2$;

2. $-\alpha_q - \frac{1}{2}(c_q - s_q - 1) > 0$ for all $q = 2, \ldots, k$ and $-\alpha_1 - \frac{1}{2}(c_1 - s_2 - 1) > 0$;
3. $-\alpha_1 - \frac{1}{2}(c_1 - s_2 + 1) - \gamma_2 > \frac{s_2 - 1}{2}$ where $\gamma_2 = \sum_{j \in J(P, S_2)}(\alpha_j - \beta_2 + \frac{c_j - s_2}{2})$.

THEOREM 3.4. *Let $G$ be a noncomplete decomposable graph and let $P$ be a perfect order of its cliques. Then $B_P \subset \mathcal{B}$. More specifically for $\theta \in Q_G$ and $(\alpha, \beta) \in B_P$ the integral (3.2) converges and*

$$\int_{P_G} e^{-\langle \theta, y \rangle} H_G(\alpha, \beta; \varphi(y)) \nu_G(dy) = \Gamma_{\text{II}}(\alpha, \beta) H_G(\alpha, \beta; \theta), \tag{3.22}$$

*where*

$$\Gamma_{\text{II}}(\alpha, \beta) = \Gamma_{s_2}\left[-\alpha_1 - \frac{c_1 - s_2}{2} - \gamma_2\right] \frac{\Gamma_{c_1}(-\alpha_1)}{\Gamma_{s_2}(-\alpha_1 - (c_1 - s_2)/2)}$$
$$\times \prod_{j=2}^{k} \frac{\Gamma_{c_j}(-\alpha_j)}{\Gamma_{s_j}(-\alpha_j - (c_j - s_j)/2)}. \tag{3.23}$$

It is interesting to reexpress (3.22) in a slightly different way. Writing $y = \hat{x}^{-1}$ with $x \in Q_G$ and recalling that the image of $\nu_G(dy)$ under $y \mapsto \varphi(y) = x$ is $\mu_G(dx)$, we see that (3.22) can be rewritten as

$$\int_{Q_G} e^{-\langle \theta, \hat{x}^{-1} \rangle} H_G(\alpha, \beta; x) \mu_G(dx) = \Gamma_{\text{II}}(\alpha, \beta) H_G(\alpha, \beta; \theta). \tag{3.24}$$

From the two theorems above, it follows immediately that

$$\mathcal{A} \supset \bigcup_P A_P \quad \text{and} \quad \mathcal{B} \supset \bigcup_P B_P,$$

where the union of all $A_P$ and all $B_P$ is taken over all possible perfect orders of the cliques of $G$. Before making some important remarks, let us give an example.

EXAMPLE. Consider the graph $G = A_4 : \overset{1}{\bullet} - \overset{2}{\bullet} - \overset{3}{\bullet} - \overset{4}{\bullet}$.

Let $C_1 = \{1, 2\}$, $C_2 = \{2, 3\}$, $C_3 = \{3, 4\}$, $S_2 = \{2\}$, $S_3 = \{3\}$, and let $\alpha(C_i) = \alpha_i$, $i = 1, 2, 3$, $\beta(S_i) = \beta_i, i = 2, 3$. Then $P_1 = (C_1, C_2, C_3)$ and $P_2 = (C_2, C_1, C_3)$ are perfect orders of the cliques. The orders $P'_1 = (C_3, C_2, C_1)$ and $P'_2 = (C_2, C_3, C_1)$ are also perfect orders analogous respectively to $P_1$ and $P_2$. On the other hand, the only other possible order $(C_1, C_3, C_2)$ and its analog $(C_3, C_1, C_2)$ are not perfect. Let us therefore identify $A_P$ and $B_P$ for $P_1, P'_1$ and $P_2, P'_2$:

$$A_{P_1} = \{(\alpha_1, \alpha_2, \alpha_3, \beta_2, \beta_3) | \alpha_i > \tfrac{1}{2}, \ i = 1, 2, 3, \ \alpha_1 + \alpha_2 - \beta_2 > 0, \ \alpha_3 = \beta_3\},$$
$$A_{P'_1} = \{(\alpha_1, \alpha_2, \alpha_3, \beta_2, \beta_3) | \alpha_i > \tfrac{1}{2}, \ i = 1, 2, 3, \ \alpha_2 + \alpha_3 - \beta_3 > 0, \ \alpha_1 = \beta_2\},$$

while $A_{P_2} = A_{P_1}$ and $A_{P'_2} = A_{P'_1}$. In a parallel way, we have

$$B_{P_1} = \{(\alpha_1, \alpha_2, \alpha_3, \beta_2, \beta_3) | -\alpha_i > 0, \ i = 1, 2, 3,$$
$$-\alpha_1 - \alpha_2 + \beta_2 - 1 > 0, \ \alpha_3 + \tfrac{1}{2} = \beta_3\},$$
$$B_{P'_1} = \{(\alpha_1, \alpha_2, \alpha_3, \beta_2, \beta_3) | -\alpha_i > 0, \ i = 1, 2, 3,$$
$$-\alpha_2 - \alpha_3 + \beta_3 - 1 > 0, \ \alpha_1 + \tfrac{1}{2} = \beta_2\},$$

while $B_{P_2} = B_{P_1}$ and $B_{P'_2} = B_{P'_1}$.

REMARKS. 1. The domains $A_P$ and $B_P$ on which (3.19) and (3.22), respectively, or equivalently (3.21) and (3.24) hold, depend upon the chosen perfect order $P$ of the cliques. Since the functions $H_G$ do not depend upon $P$, it is clear that, even though the expressions of $\Gamma_{\text{I}}$ and $\Gamma_{\text{II}}$ depend upon $P$, their values do not.

2. Since assumption 1 of Theorems 3.3 and 3.4 represents $k' - 1$ constraints on the set of $(\alpha, \beta)$'s, we see that in general each set $A_P$ is of dimension $k + 1$.

3. From Theorems 3.3 and 3.4, the integrals (3.1) and (3.2) are finite and constant multiples of $H_G(\alpha, \beta; \sigma)$ and $H_G(\alpha, \beta; \theta)$ for $(\alpha, \beta)$ in $\bigcup_P A_P$ and $\bigcup_P B_P$, respectively. Using Hölder's inequality it is immediate to prove that these integrals are also finite on the convex hull of $\bigcup_P A_P$ and $\bigcup_P B_P$. So the question naturally arises as to whether $\mathcal{A}$ and $\mathcal{B}$ are larger than $\bigcup_P A_P$ and $\bigcup_P B_P$. We only have a partial answer to this. We have seen in the previous section that, when $G$ is homogeneous, $\mathcal{A}$ and $\mathcal{B}$ are completely known and of full dimension $k + k'$. However, if we consider the homogeneous example given in Section 3.3 and treat it using the methods given in this section, we will find that the 24 possible orders are all perfect with $P_1 = (C_1, C_2, C_3, C_4)$ and $P_2 = (C_1, C_4, C_2, C_3)$ being the only perfect orders yielding distinct $A_P$'s. We have

$$A_{P_1} = \{(\alpha, \beta) | \alpha_i > 1, \ i = 1, 2, 3, \ \alpha_4 > \tfrac{1}{2}, \ \alpha_1 + \alpha_2 + \alpha_3 - 2\beta_2 > \tfrac{1}{2}, \ \alpha_4 = \beta_1\},$$
$$B_{P_1} = \{(\alpha, \beta) | \alpha_i < 0, \ i = 1, 2, 3, 4,$$
$$-\alpha_1 - \alpha_2 - \alpha_3 + 2\beta_2 > \tfrac{5}{2}, \ \alpha_4 - \beta_1 + \tfrac{1}{2} = 0\},$$
$$A_{P_2} = \{(\alpha, \beta) | \alpha_i > 1, \ i = 1, 2, 3,$$
$$\alpha_4 > \tfrac{1}{2}, \alpha_1 + \alpha_4 - \beta_1 > 0, \ \alpha_2 + \alpha_3 - 2\beta_2 = 0\},$$
$$B_{P_2} = \{(\alpha, \beta) | \alpha_i < 0, \ i = 1, 2, 3, 4,$$
$$-\alpha_1 - \alpha_4 + \beta_1 > 2, \ \alpha_2 + \alpha_3 - 2\beta_2 + 1 = 0\}.$$

Clearly, $A_{P_1} \cup A_{P_2}$ is included in, but not equal to, $\mathcal{A}$ as given in (3.17). Similarly, $B_{P_1} \cup B_{P_2}$ is included in, but not equal to, $\mathcal{B}$ as given in (3.18). The question is therefore whether in the nonhomogeneous case it is possible to identify $\mathcal{A}$ and $\mathcal{B}$. In the next section, we find $\mathcal{A}$ and $\mathcal{B}$ for $G = A_4$ and we see that they are of dimension strictly less than $k + k'$. Thus $G = A_4$ is a counterexample to the hypothesis that in the nonhomogeneous case we could also define a set of dimension $k + k'$ on which (3.21) and (3.22) hold.

### 3.5. The case $G = A_4$.

Let $G$ be $A_4$ as in the previous example. Then we write

$$\sigma = \begin{pmatrix} \sigma_1 & \sigma_{12} & & \\ \sigma_{21} & \sigma_2 & \sigma_{23} & \\ & \sigma_{32} & \sigma_3 & \sigma_{34} \\ & & \sigma_{43} & \sigma_4 \end{pmatrix}$$

for $\sigma \in Q_G$, with $\sigma_{ij} = \sigma_{ji}$, $\sigma_{i.j} = \sigma_i - \sigma_{ij}\sigma_j^{-1}\sigma_{ji}$ and similarly for $\theta \in Q_G$.

PROPOSITION 3.2. *Consider the graph $G = A_4$ with cliques and separators*

$$C_1 = \{1,2\}, \quad C_2 = \{2,3\}, \quad C_3 = \{3,4\}, \quad S_2 = \{2\}, \quad S_3 = \{3\}.$$

*Let $\alpha_i = \alpha(C_i)$, $i = 1,2,3$, $\beta_i = \beta(S_i)$, $i = 2,3$. Define*

$$\mathcal{A}_4 = \{(\alpha,\beta) | \alpha_i > \tfrac{1}{2}, \ i = 1,2,3, \ \alpha_1 + \alpha_2 > \beta_2, \ \alpha_2 + \alpha_3 > \beta_3\}.$$

*Then the following integral converges for all $\sigma \in Q_{A_4}$ if and only if $(\alpha,\beta)$ is in $\mathcal{A}_4$. Under these conditions, it is equal to*

$$\int_{Q_G} e^{-\langle x, \hat{\sigma}^{-1} \rangle} H_G(\alpha, \beta; x) \mu_G(dx)$$

$$= \pi^{3/2} \Gamma\left(\alpha_1 - \frac{1}{2}\right) \Gamma\left(\alpha_2 - \frac{1}{2}\right) \Gamma\left(\alpha_3 - \frac{1}{2}\right) \Gamma(\alpha_1 + \alpha_2 - \beta_2)$$

(3.25)
$$\times \Gamma(\alpha_2 + \alpha_3 - \beta_3)(\Gamma(\alpha_2))^{-1} \sigma_{1\cdot 2}^{\alpha_1} \sigma_{2\cdot 3}^{\alpha_1 + \alpha_2 - \beta_2} \sigma_{3\cdot 2}^{\alpha_2 + \alpha_3 - \beta_3} \sigma_{4\cdot 3}^{\alpha_3}$$

$$\times {}_2F_1\left(\alpha_1 + \alpha_2 - \beta_2, \alpha_2 + \alpha_3 - \beta_3, \alpha_2, \frac{\sigma_{23}^2}{\sigma_2 \sigma_3}\right),$$

*where ${}_2F_1$ denotes the hypergeometric function. Similarly we define*

$$\mathcal{B}_4 = \{(\alpha,\beta) | \alpha_1 < 0, \ \alpha_3 < 0, \ \beta_2 - \alpha_1 - \alpha_2 - \tfrac{1}{2} > 0,$$
$$\beta_3 - \alpha_2 - \alpha_3 - \tfrac{1}{2} > 0, \ \beta_2 + \beta_3 - \alpha_1 - \alpha_2 - \alpha_3 - \tfrac{3}{2} > 0\}.$$

*Then the following integral converges for all for $\theta \in Q_G$ if and only if $(\alpha,\beta) \in \mathcal{B}_4$. Under these conditions, it is equal to*

$$\int_{P_G} e^{-\langle \theta, y \rangle} H_G(\alpha, \beta; \varphi(y)) \nu_G(dy)$$
$$= \pi^{3/2} \theta_{1\cdot 2}^{\alpha_1} \theta_2^{\alpha_1+\alpha_2-\beta_2} \theta_3^{\alpha_2+\alpha_3-\beta_3} \theta_{4\cdot 3}^{\alpha_3}$$
$$\times \Gamma(-\alpha_1) \Gamma\left(\beta_2 - \alpha_1 - \alpha_2 - \frac{1}{2}\right) \Gamma\left(\beta_2 + \beta_3 - \alpha_1 - \alpha_2 - \alpha_3 - \frac{3}{2}\right)$$
(3.26) $\quad \times \Gamma\left(\beta_3 - \alpha_2 - \alpha_3 - \frac{1}{2}\right) \Gamma(-\alpha_3) (\Gamma(\beta_2 + \beta_3 - \alpha_1 - \alpha_2 - \alpha_3 - 1))^{-1}$
$$\times {}_2F_1\left(\beta_2 - \alpha_1 - \alpha_2 - \frac{1}{2}, \beta_3 - \alpha_2 - \alpha_3 - \frac{1}{2};\right.$$
$$\left.\beta_2 + \beta_3 - \alpha_1 - \alpha_2 - \alpha_3 - 1; \frac{\theta_{23}^2}{\theta_2 \theta_3}\right).$$

The results above are obtained by a nontrivial and long computation. A central part of this computation is the following lemma.

LEMMA 3.1. *Consider the following $2 \times 2$ random matrix $X = \begin{bmatrix} X_1 & X_{12} \\ X_{12} & X_2 \end{bmatrix}$ with the Wishart distribution*
$$w_2(p, c^{-1}; dx) = \frac{(\det c)^p}{\Gamma_2(p)} e^{-\langle x, c \rangle} (x_1 x_2 - x_{12}^2)^{p-3/2} \mathbf{1}_{M_2^+}(x) \, dx_1 \, dx_2 \, dx_{12}$$
*with $p \geq 1/2$ and $c = \begin{bmatrix} c_1 & c_{12} \\ c_{12} & c_2 \end{bmatrix}$ positive definite. For $a_1 > -p$ and $a_2 > -p$, the Mellin transform of $(X_1, X_2)$ is*
$$E(X_1^{a_1} X_2^{a_2}) = \frac{(\det c)^p}{c_1^{a_1+p} c_2^{a_2+p}} \frac{\Gamma(a_1+p)\Gamma(a_2+p)}{\Gamma(p)^2} {}_2F_1\left(a_1+p, a_2+p; p; \frac{c_{12}^2}{c_1 c_2}\right).$$

The proofs of Proposition 3.2 and Lemma 3.1 are omitted.

We now derive from Proposition 3.2 the sets $\mathcal{A}$ and $\mathcal{B}$ when $G = A_4$.

COROLLARY 3.1. *Let $G = A_4$. Then $\mathcal{A} = \bigcup A_P$ and $\mathcal{B} = \bigcup B_P$, where the unions are taken over the two possible $A_P$ and $B_P$. The dimension of $\mathcal{A}$ and $\mathcal{B}$ is therefore strictly less than $k + k'$.*

PROOF. Since the two statements are quite similar, we prove the second one only. We use the equality (see [1], Formula 15.3.3)
$$(3.27) \qquad (1-z)^{a+b-c} {}_2F_1(a, b; c; z) = {}_2F_1(c-a, c-b; c; z).$$
Using (3.27) on the right hand side of (3.26) above with $z = \frac{\theta_{23}^2}{\theta_2 \theta_3}$ and
$$(a, b, c) = (\beta_2 - \alpha_1 - \alpha_2 - \tfrac{1}{2}, \beta_3 - \alpha_2 - \alpha_3 - \tfrac{1}{2}; \beta_2 + \beta_3 - \alpha_1 - \alpha_2 - \alpha_3 - 1),$$

we see that

$$\frac{1}{H_G(\alpha,\beta;\theta)}\int_{Q_G} e^{-\langle\theta,y\rangle} H_G(\alpha,\beta;\varphi(y))\nu_G(dy) = C\,{}_2F_1(c-a,c-b;c;z),$$

where $C$ is a constant which does not depend on $z$. Now clearly from its Taylor expansion the hypergeometric function $z \mapsto {}_2F_1(c-a,c-b;c;z)$ is a constant if and only if either $c-a$ or $c-b$ is zero, which together with $\mathcal{B}_4$ proves that constancy occurs if and only if $(\alpha,\beta)$ belongs to one of the two possible $B_P$ as given in the example of Section 3.4 above. □

REMARK. One can prove that the convex hulls in $\mathbb{R}^5$ of $\mathcal{A}$ and $\mathcal{B}$ are, respectively, strictly included in $\mathcal{A}_4$ and $\mathcal{B}_4$ as defined in Proposition 3.2.

**4. Properties of the Type I and II Wisharts.** Let us recall that for a given decomposable graph $G$, the $r$-dimensional graphical Gaussian model Markov with respect to $G$ is the family of distributions

$$\mathcal{N}_G = \{N_r(0,\Sigma), \Sigma \in Q_G\}.$$

Dawid and Lauritzen [7], page 1306, have shown that this model is strong meta Markov. This can also be shown directly since, using the notation of (2.8), for a given perfect order of the cliques and with the convention that $x_{[1]\cdot} = x_{C_1}, \Sigma_{[1]\cdot} = \Sigma_{C_1}, x_{<1>} = 0, r_1 = c_1$, the density of $X \in \mathbb{R}^r$ with distribution $N_r(0,\Sigma) \in \mathcal{N}_G$ can be written as

$$(4.1)\ f(x) = \prod_{i=1}^{k}\frac{1}{(2\pi)^{r_i/2}\Sigma_{[i]\cdot}^{1/2}} e^{(-1/2)\langle(x_{[i]}-x_{[i>}\Sigma_{<i>}^{-1}\Sigma_{<i]}),\Sigma_{[i]\cdot}^{-1}(x_{[i]}-x_{[i>}\Sigma_{<i>}^{-1}\Sigma_{<i]})\rangle}.$$

The parameters

$$(4.2)\qquad F_1 = \Sigma_{C_1} \qquad (L_i = \Sigma_{[i>}\Sigma_{<i>}^{-1}, N_i = \Sigma_{[i]\cdot}),\qquad i=2,\ldots,k,$$

of the distributions of $X_{C_1}, X_{[i]}|x_{H_{i-1}}, i=2,\ldots,k$, respectively, are clearly variation independent in the sense that any parameter $(\Sigma_{[i>}\Sigma_{<i>}^{-1},\Sigma_{[i]\cdot})$ of the distribution of $X_{[i]}|X_{H_{i-1}}$ is compatible with any parameter $\{F_1,(L_j,N_j),j<i\}$ of $X_{H_{i-1}}$ and any parameter $\{L_j,N_j,j>i\}$ of $X|X_{H_i}$. It then follows that for a decomposition $(A,B)$ of $G$ the parameter $\Sigma_A$ of the distribution of $X_A$ is variation independent of the parameter $\Sigma_{B|A}$ of the conditional distribution of $X_B$ given $X_A$. This, according to Definition 4.3 of [7], means that the model $\mathcal{N}_G$ is strong meta-Markov. For this model, Dawid and Lauritzen ([7], page 1306 and page 1308) have shown that the distribution of the maximum likelihood estimator of $\Sigma$, that is the hyper Wishart, is weak hyper Markov and that the hyper inverse Wishart, the inverse of the $G$-Wishart, is a conjugate prior on $\Sigma$ which is strong hyper Markov.

We are now going to show parallel results for the Type I and II Wisharts: the Type I Wishart is weak hyper Markov, the inverse of the Type II Wishart forms a conjugate family for the scale parameter of the $\mathcal{N}_G$ model and for any direction given to the graph by a perfect order of its cliques, and the inverse Type II Wishart is strong directed hyper Markov. Since we have seen in Section 3.2 that the hyper Wishart is a particular case of the type I Wishart and the hyper inverse Wishart is a particular case of the inverse Type II, it is not surprising that their generalizations hold parallel Markov properties.

One might wonder whether the term "hyper" is adequate when talking about the weak Markov property of the Type I Wishart since this distribution has so far neither been identified as the distribution of an estimator nor as a prior distribution for the parameter of a Gaussian model. It is certainly adequate for the inverse of the Type II Wishart since, as we are going to prove right away in Section 4.1, it forms a conjugate family of prior distributions for the scale parameter of the $\mathcal{N}_G$ model, with a shape parameter set of dimension at least $k+1$. We will then prove the hyper Markov properties in Section 4.2. Our main results below are Corollary 4.1 and Theorem 4.4.

4.1. *Conjugate prior distributions.* The family of inverse Type II Wishart distributions has several properties that make it useful as a rich family of conjugate prior distributions for the scale parameter $\Sigma_G$ of the graphical Gaussian model Markov with respect to a decomposable graph $G$. Recall that, following the notation used in the Introduction, if $\Sigma$ is the positive definite covariance matrix for $N_r(0, \Sigma) \in \mathcal{N}_G$, then $\Sigma_G = \pi(\Sigma) \in Q_G$ is the scale parameter for the $N_r(0, \Sigma)$ distribution. We have the following general result.

THEOREM 4.1. *Let $G$ be a decomposable graph and let $P$ be a perfect order of its cliques. Let $D$ be in $Q_G$, let $(\alpha, \beta)$ be in $A_P$ and $(\alpha', \beta')$ be in $B_P$. If the joint distribution of $(X, \Sigma_G)$ on $Q_G \times Q_G$ is $W_{Q_G}(\alpha, \beta, \sigma; dx)IW_{P_G}(\alpha', \beta', D; d\sigma)$, then the conditional distribution of $\Sigma_G$ knowing $X = x$ is $IW_{P_G}(\alpha' - \alpha, \beta' - \beta, D + x; d\sigma)$ and the marginal distribution of $X$ is an $F$ distribution of the first kind with parameter $(\alpha, \beta, \alpha', \beta', D)$.*

PROOF. The joint distribution of $(X, \Sigma_G)$ is

$$W_{Q_G}(\alpha, \beta, \sigma; dx)IW_{P_G}(\alpha', \beta', D; d\sigma)$$

$$= \frac{e^{-\langle x, \hat{\sigma}^{-1}\rangle} H_G(\alpha, \beta; x)}{\Gamma_{\mathrm{I}}(\alpha, \beta)H_G(\alpha, \beta; \sigma)}\mu_G(dx)\frac{e^{-\langle D, \hat{\sigma}^{-1}\rangle} H_G(\alpha', \beta'; \sigma)}{\Gamma_{\mathrm{II}}(\alpha', \beta')H_G(\alpha', \beta'; D)}\mu_G(d\sigma)$$

$$= \left[\frac{H_G(\alpha, \beta; x)}{\Gamma_{\mathrm{I}}(\alpha, \beta)\Gamma_{\mathrm{II}}(\alpha', \beta')H_G(\alpha', \beta'; D)}\mu_G(dx)\right]$$

$$\times e^{-\langle x+D, \hat{\sigma}^{-1}\rangle} H_G(\alpha' - \alpha, \beta' - \beta; \sigma)\mu_G(d\sigma),$$

from which the result follows immediately. □

Theorem 4.1 shows that the family of $IW_{P_G}$ distributions is a conjugate family for the scale parameter $\sigma$ of the $W_{Q_G}(\alpha, \beta, \sigma; dx)$. Consider now a sample $Z_1, \ldots, Z_n$ from a Gaussian distribution Markov with respect to $G$, and let $S = \frac{1}{n} \sum_{i=1}^n Z_i Z_i^t$. Then $\pi(S)$, the maximum likelihood estimator of $\Sigma_G$ is such that $n\pi(S)$ is hyper Wishart with shape parameter $p = \frac{n}{2}$ and scale parameter $\Sigma_G$, that is, Wishart of Type I with shape parameter

$$\alpha(C) = p, \qquad C \in \mathcal{C}, \qquad \beta(S) = p, \qquad S \in \mathcal{S}, \qquad p > \max_{C \in \mathcal{C}} \frac{|C|-1}{2}$$

and scale parameter $2\Sigma_G$. Applying Theorem 4.1 to $X = n\pi(S)$, we obtain the following corollary.

COROLLARY 4.1. *Let $G$ be decomposable and let $P$ be a perfect order of its cliques. Let $(Z_1, \ldots, Z_n)$ be a sample from the $N_r(0, \Sigma)$ distribution with $\Sigma_G \in Q_G$. If the prior distribution on $2\Sigma_G$ is $IW_{P_G}(\alpha', \beta', D)$ with $(\alpha', \beta') \in B_P$ and $D \in Q_G$, the posterior distribution of $2\Sigma_G$, given $nS = \sum_{i=1}^n Z_i Z_i^t$, is $IW_{P_G}(\alpha' - \frac{n}{2}, \beta' - \frac{n}{2}, D + \pi(nS))$, where $\alpha' - \frac{n}{2} = (\alpha_1 - \frac{n}{2}, \ldots, \alpha_k - \frac{n}{2})$ and $\beta' - \frac{n}{2} = (\beta'_1 - \frac{n}{2}, \ldots, \beta'_k - \frac{n}{2})$ are such that $(\alpha' - \frac{n}{2}, \beta' - \frac{n}{2}) \in B_P$ and $D + \pi(nS) \in Q_G$.*

This means that the family $\{IW_{P_G}(\alpha, \beta, D), (\alpha, \beta) \in \mathcal{A}, D \in Q_G\}$ is a conjugate family for the scale parameter $\Sigma_G$ of the Gaussian model Markov with respect to $G$. We note this family has its shape parameter set of dimension at least $k + 1$ and is therefore much richer than the traditional Diaconis–Ylvisaker family with shape parameter set of dimension equal to 1.

Theorem 4.1 can also be immediately transcribed to the homogeneous case using the variables $(x_{[u]}, x_{[u>} x_{<u>}^{-1}, u \in T)$ and we obtain the following result.

THEOREM 4.2. *Let $G$ be a homogeneous graph. Let $D$ be in $Q_G$, let $(\alpha, \beta)$ be in $\mathcal{A}$ and $(\alpha', \beta')$ be in $\mathcal{B}$. If $(X, \Sigma_G) \in Q_G \times Q_G$ and the joint distribution of $(X_{[u]}, X_{[u>} X_{<u>}^{-1}, \Sigma_{[u]}, \Sigma_{[u>} \Sigma_{<u>}^{-1}, u \in T)$ is $W^*_{Q_G}(\alpha, \beta, \sigma) IW^*_{P_G}(\alpha', \beta', D)$, then the conditional distribution of $\Sigma$ knowing $X = x$ is $IW^*_{P_G}(\alpha' - \alpha, \beta' - \beta, D + x)$.*

We now have the following result dual to Theorem 4.1.

THEOREM 4.3. *Let $P$ be a perfect order of the cliques of $G$. Let $\sigma$ be in $Q_G$. Let $(\alpha, \beta)$ be in $B_P$ and $(\alpha', \beta')$ be in $A_P$. If the joint distribution of $(Y, \Theta)$ on $P_G \times P_G$ is $W_{P_G}(\alpha, \beta, \theta) IW_{Q_G}(\alpha', \beta', \sigma)$, then the conditional distribution of $\Theta$ knowing $Y = y$ is $IW_{Q_G}(\alpha' - \alpha, \beta' - \beta, \varphi(y + \hat{\sigma}^{-1}))$.*

Theorem 4.3 shows that the family of $IW_{Q_G}$ distributions is a conjugate family for the scale parameter $\theta$ of the $W_{P_G}(\alpha, \beta, \theta)$.

4.2. *Markov properties.* We now want to show that the Type II inverse Wishart $IW_{P_G}$ is strong directed hyper Markov and the Type I Wishart weak hyper Markov. Let $M(G)$ denote the set of all Markov probabilities over $G$. Let the distribution $P_\rho \in M(G)$ be parametrized by $\rho$. Now we randomize $\rho$ according to a law $\mathcal{L}(\rho)$. For any subset $A$ of $V$, let $\rho_A$ denote the parameter of the marginal distribution of $X_A$ (more specifically, we should write $\rho \sim_A \rho'$ if the marginal distributions of $X_A$ under $P_\rho$ and $P_{\rho'}$ coincide and call $\rho_A$ the equivalence class of $\rho$ for the equivalence relation $\sim_A$). The parameter $\rho_{A|B}$ of the conditional distribution of $X_A$ knowing $X_B$ could be defined in a similar way. We say that $\mathcal{L}(\rho)$ is weak hyper Markov over $G$ if under $\mathcal{L}(\rho)$, for any decomposition $(A, B)$ of $V$,

$$\rho_A \perp\!\!\!\perp \rho_B | \rho_{A \cap B}. \tag{4.3}$$

We say that $\mathcal{L}(\rho)$ is strong hyper Markov over $G$ if, under $\mathcal{L}(\rho)$, for any decomposition $(A, B)$ of $V$,

$$\rho_{A|B} \perp\!\!\!\perp \rho_B. \tag{4.4}$$

Let $P$ be any perfect order of the cliques and consider a perfect numbering of the vertices compatible with $P$ (see [11], page 18). Let $D$ be the directed graph obtained from $G$ by directing all edges in $G$ from the vertex with the smallest number to the vertex with the largest number. We say that a law $\mathcal{L}(\rho)$ is weak directed hyper Markov over $D$ if for all $v \in V$,

$$\rho_v \perp\!\!\!\perp \rho_{pr(v)} | \rho_{pa(v)}, \tag{4.5}$$

where $pa(v)$ denotes the sets of parents of $v$ in $D$ and $pr(v)$ denotes the sets of predecessors, that is, the vertices with a lower number than $v$.

We say that $\mathcal{L}(\rho)$ is strong directed hyper Markov over $D$ if for all $v \in V$,

$$\rho_{v|pa(v)} \perp\!\!\!\perp \rho_{pr(v)}. \tag{4.6}$$

Let us also recall that a random variable on $M_r^+$ is said to follow the Wishart $w_r(p, \sigma)$ distribution if its density with respect to the Lebesgue measure is

$$\frac{|x|^{p-(r+1)/2}}{|\sigma|^p \Gamma_r(p)} e^{-\langle x, \sigma^{-1} \rangle},$$

in which case its inverse $U = X^{-1}$ is said to follow the inverse Wishart distribution $iw_r(p, \theta)$, where $\theta = \sigma^{-1}$, with density with respect to the Lebesgue measure equal to

$$\frac{|u|^{-p-(r+1)/2}}{|\sigma|^p \Gamma_r(p)} e^{-\langle u^{-1}, \theta \rangle}.$$

Finally, we will use the notation $x_{[12>}$ and $x_{[1]\cdot}$ for

$$x_{[12>} = x_{C_1 \setminus S_2, S_2} x_{S_2}^{-1} \quad \text{and} \quad x_{[1]\cdot} = x_{C_1 \setminus S_2 \cdot S_2}.$$

THEOREM 4.4. *Let $G$ be a decomposable graph $G$ and let $P$ be a perfect order of its cliques. Then, for $(\alpha, \beta) \in B_P$ and for the direction given by $P$, the inverse Type* II *Wishart is strong directed hyper Markov. More precisely, if $X \sim IW_{P_G}(\alpha, \beta, \theta)$ with $(\alpha, \beta) \in B_P$ and $\theta \in Q_G$, then, with the convention that $s_1 = s_2$,*

$$x_{[i]\cdot} \sim iw_{c_i - s_i}(-\alpha_i, \theta_{[i]\cdot}), \qquad i = 1, \ldots, k,$$

$$x_{[12>} | x_{[1]\cdot} \sim N_{(c_1 - s_2) \times s_2}(\theta_{[12>}, 2\theta_{<2>}^{-1} \otimes x_{[1]\cdot}),$$

$$x_{<2>} \sim iw_{s_2}\left(-\left(\alpha_1 + \frac{c_1 - s_2}{2} + \gamma_2\right), \theta_{<2>}\right),$$

$$x_{[j>} x_{<j>}^{-1} | x_{[j]\cdot} \sim N_{(c_j - s_j) \times s_j}(\theta_{[j>} \theta_{<j>}^{-1}, 2\theta_{<j>}^{-1} \otimes x_{[j]\cdot}), \qquad j = 2, \ldots, k$$

*and*

(4.7) $$\{(x_{[12>}, x_{[1]\cdot}), x_{<2>}, (x_{[j>} x_{<j>}^{-1}, x_{[j]\cdot}), \ j = 2, \ldots, k\}$$

*are mutually independent.*

PROOF. From (4.1), we know that $(\Sigma_{[i>} \Sigma_{<i>}^{-1}, \Sigma_{[i]\cdot})$ is the parameter of the distribution of $Z_{[i]}$ given $Z_{<i>}$ when $Z \sim N_r(0, \Sigma) \in \mathcal{N}_G$. Therefore it follows from the remark following Theorem 2.6 and Proposition 3.8 of [7] that to construct a weak hyper Markov distribution for $\rho = \Sigma_G \in Q_G$, it is sufficient to build a weak directed hyper Markov distribution for a given direction of the vertices compatible with a given perfect order of the cliques, that is, a distribution with density of the form

$$p(\Sigma_G) = p_{C_1}(\Sigma_{C_1}) \prod_{i=2}^{k} p_i(\Sigma_{[i>} \Sigma_{<i>}^{-1}, \Sigma_{[i]\cdot} | \Sigma_{<j>})$$

(4.8) $$= p_1(\Sigma_{[12>}, \Sigma_{[1]\cdot} | \Sigma_{<2>}) p_{<2>}(\Sigma_{<2>})$$

$$\times \prod_{i=2}^{k} p_i(\Sigma_{[i>} \Sigma_{<i>}^{-1}, \Sigma_{[i]\cdot} | \Sigma_{<j>}).$$

If we want to show that, for the given direction, this distribution is in fact strong directed hyper Markov, by Proposition 3.13 of [7] it is sufficient to show that

(4.9) $$\{(\Sigma_{[12>}, \Sigma_{[1]\cdot}), \ \Sigma_{<2>}, \ (\Sigma_{[j>} \Sigma_{<j>}^{-1}, \Sigma_{[j]\cdot}), \ j = 2, \ldots, k\}$$

are mutually independent. Let us now show that the Inverse Type II Wishart satisfies both (4.8) and (4.9). Let $X \sim IW_{P_G}(\alpha, \beta, \theta)$ with $(\alpha, \beta) \in B_P$ and $\theta \in Q_G$. Combining (A.9) and (A.11) of the Appendix, we show that the image of the $IW_{P_G}(\alpha, \beta, \theta; dx)$ distribution by the change of variables (2.12) and (A.5) is

$$IW_{P_G}^{**}(\alpha, \beta, \theta; dx_{[1]\cdot}, dx_{[12>}, d(x_{S_2}), d(x_{[j>}x_{<j>}^{-1}), dx_{[j]\cdot}, j = 2, \ldots, k)$$

(4.10) $$\propto e^{-\langle x_{[1]\cdot}^{-1}, \theta_{[1]\cdot}\rangle}|x_{[1]\cdot}|^{\alpha_1 - (c_1 - s_2 + 1)/2}|x_{[1]\cdot}|^{-s_2/2}$$

$$\times e^{-\langle (x_{[12>} - \theta_{[12>}), x_{[1]\cdot}^{-1}(x_{[12>} - \theta_{[12>})\theta_{<2>}\rangle}$$

(4.11) $$\times e^{-\langle x_{<2>}^{-1}, \theta_{<2>}\rangle}|x_{<2>}|^{\alpha_1 + (c_1 - s_2)/2 + \gamma_2 - (s_2 + 1)/2}\, dx_{<2>}$$

(4.12) $$\times \prod_{j=2}^{k} |x_{[j]\cdot}|^{\alpha_j - (c_j - s_j + 1)/2} e^{-\langle x_{[j]\cdot}^{-1}, \theta_{[j]\cdot}\rangle}|x_{[j]\cdot}|^{-s_j/2}$$

(4.13) $$\times \prod_{j=2}^{k} e^{-\langle (x_{[j>}x_{<j>}^{-1} - \theta_{[j>}\theta_{<j>}^{-1}), x_{[j]\cdot}^{-1}(x_{[j>}x_{<j>}^{-1} - \theta_{[j>}\theta_{<j>}^{-1})\theta_{<j>}\rangle}$$

$$\times dx_{[12>}\, dx_{[1]\cdot} \prod_{j=2}^{k} d(x_{[j>}x_{<j>}^{-1})\, dx_{[j]\cdot}.$$

We see that the densities (4.10), (4.11), (4.12) and (4.13) with $x$ replaced by $\Sigma_G$ are exactly of the form required for the respective factors of (4.8). It follows that the $IW_{P_G}(\alpha, \beta, \theta)$ is weak directed hyper Markov, but we also see from (4.10)–(4.13) above that the independence in (4.9) is satisfied and therefore $IW_{P_G}(\alpha, \beta, \theta)$ is strong directed hyper Markov. The densities of

$$(x_{[1]\cdot}, x_{[12>}, x_{<2>}, x_{[j]\cdot}, x_{[j>}x_{<j>}^{-1}, j = 2, \ldots, k)$$

are also clearly as indicated in the theorem. □

This strong directed hyper Markov property of the Type II inverse Wishart corresponds to the strong hyper Markov property for the inverse $G$-Wishart, that is, the hyper inverse Wishart. We do not quite have the strong hyper Markov property because the parameters $(\alpha, \beta) \in B_P$ are linked to the perfect order $P$. The property analogous to the weak hyper Markov property of the hyper Wishart is given in the following theorem.

THEOREM 4.5. *Let $G$ be a decomposable graph and let $P$ be a perfect order of its cliques. Then, for $(\alpha, \beta) \in A_P$, the Type I Wishart is weak hyper Markov. More precisely, if $X \sim W_{Q_G}(\alpha, \beta, \sigma)$ with $(\alpha, \beta) \in A_P$ and $\sigma \in Q_G$, then*

$$x_{[1]\cdot} \sim w_{c_1 - s_2}\left(\alpha_1 - \frac{s_2}{2}, \sigma_{[1]\cdot}\right),$$

$$x_{[12>}|x_{<2>} \sim N_{(c_1-s_2) \times s_2}(\sigma_{[12>}, 2x_{<2>}^{-1} \otimes \sigma_{[1].}),$$

$$x_{<2>} \sim w_{s_2}(\alpha_1 + \delta_2, \sigma_{<2>}),$$

$$x_{[j>}x_{<j>}^{-1}|x_{<j>} \sim N_{(c_j-s_j) \times s_j}(\sigma_{[j>}\sigma_{<j>}^{-1}, 2x_{<j>}^{-1} \otimes \sigma_{[j].}),$$

$$x_{[j].} \sim w_{c_j-s_j}\left(\alpha_j - \frac{s_j}{2}, \sigma_{[j].}\right), \qquad j = 2, \ldots, k.$$

PROOF. Using (A.4) and (A.6) of the Appendix, we see that the image of the $W_{Q_G}(\alpha, \beta, \sigma; dx)$ distribution by the change of variables (2.12) and (A.5) is

$$W_{Q_G}^{**}(\alpha, \beta, \sigma; dx_{[1].}, dx_{[12>}, dx_{<2>}, dx_{[j].}, d(x_{[j>}x_{<j>}^{-1})), j = 2, \ldots, k)$$

$$(4.14) \quad \propto |x_{[1].}|^{\alpha_1 - s_2/2 - (c_1 - s_2 + 1)/2} e^{-\langle x_{[1].}, \sigma_{[1].}^{-1} \rangle} dx_{[1].}$$

$$(4.15) \quad \times |x_{<2>}|^{(c_1-s_2)/2} e^{-\langle (x_{[12>} - \sigma_{[12>})x_{<2>}(x_{<21]} - \sigma_{<21]})\sigma_{[1].}^{-1} \rangle} dx_{[12>}$$

$$(4.16) \quad \times |x_{<2>}|^{\alpha_1 + \delta_2 - (s_2+1)/2} e^{-\langle x_{<2>}, \sigma_{<2>}^{-1} \rangle} dx_{<2>}$$

$$(4.17) \quad \times \prod_{j=2}^{k} |x_{[j].}|^{\alpha_j - s_j/2 - (c_j - s_j + 1)/2} e^{-\langle x_{[j].}, \sigma_{[j].}^{-1} \rangle} dx_{[j].}$$

$$\times |x_{<j>}|^{(c_j - s_j)/2}$$

$$(4.18) \quad \times e^{-\langle (x_{[j>}x_{<j>}^{-1} - \sigma_{[j>}\sigma_{<j>}^{-1})x_{<j>}(x_{<j>}^{-1}x_{<j]} - \sigma_{<j>}^{-1}\sigma_{<j]})\sigma_{[j].}^{-1} \rangle}$$

$$\times d(x_{[j>}x_{<j>}^{-1}).$$

From the expression of the density $W_{Q_G}^{**}(\alpha, \beta, \sigma)$ above, we see that, for $x = \Sigma_G$, (4.14) and (4.15) give $p_1(\Sigma_{[C_1 \setminus S_2, S_2>}\Sigma_{S_2}^{-1}, \Sigma_{C_1 \setminus S_2 \cdot S_2}|\Sigma_{S_2})$ of (4.8) while (4.16) gives $p_{S_2}(\Sigma_{S_2})$ and (4.17) and (4.18) give $p_i(\Sigma_{[i>}\Sigma_{<i>}^{-1}, \Sigma_{[i].}|\Sigma_{<j>})$. Therefore the $W_{Q_G}(\alpha, \beta, \sigma)$ Type I Wishart is weak directed hyper Markov and, by Proposition 3.8 of [7], weak hyper Markov. □

We note that in the proof above, the density $p_1(\Sigma_{[C_1 \setminus S_2, S_2>}\Sigma_{S_2}^{-1}, \Sigma_{C_1 \setminus S_2 \cdot S_2}|\Sigma_{S_2})$ depends upon $\Sigma_{S_2}$ and the density $p_i(\Sigma_{[i>}\Sigma_{<i>}^{-1}, \Sigma_{[i].}|\Sigma_{<j>})$ depends upon $\Sigma_{<j>}$ and therefore the $W_{Q_G}(\alpha, \beta, \sigma)$ Type I Wishart is not strong directed hyper Markov.

4.3. *Laplace transforms and expected values.* For $(\alpha, \beta) \in \mathcal{A}$, $\mathcal{F}_{(\alpha,\beta),I} = \{W_{Q_G}(\alpha, \beta, \sigma; dx), \sigma \in Q_G\}$ is the natural exponential family generated by the measure

$$\mu_{(\alpha,\beta),G}(dx) = \frac{H_G(\alpha, \beta; x)}{\Gamma_I(\alpha, \beta)} \mu_G(dx).$$

For $-y \in P_G$, the Laplace transform of $\mu_{(\alpha,\beta),G}$ is

$$L_{\mu_{(\alpha,\beta),G}}(y) = \int_{Q_G} e^{\text{tr}(xy)} \mu_{(\alpha,\beta),G}(dx) = H_G(\alpha,\beta;-\varphi(y)).$$

This is a reformulation of (3.19). It implies that, for $-y + \hat{\sigma}^{-1} \in P_G$, the Laplace transform of $W_{Q_G}(\alpha,\beta,\sigma)$ is defined by

$$(4.19) \quad \int_{Q_G} e^{\text{tr}(xy)} W_{Q_G}(\alpha,\beta,\sigma;dx) = \frac{H_G(\alpha,\beta;\varphi(\hat{\sigma}^{-1}-y))}{H_G(\alpha,\beta;\sigma)}.$$

Suppose that $(\alpha,\beta)$ and $(\alpha',\beta')$ in $\mathcal{A}$ are such that $(\alpha+\alpha',\beta+\beta')$ is still in $\mathcal{A}$. This is, for example, true for any $G$ if $(\alpha,\beta)$ and $(\alpha',\beta')$ are in the same $A_P$ and it is always true if $G$ is homogeneous. We then have the convolution formula

$$W_{Q_G}(\alpha,\beta,\sigma) * W_{Q_G}(\alpha',\beta',\sigma) = W_{Q_G}(\alpha+\alpha',\beta+\beta',\sigma),$$

a result which would be difficult to prove using densities alone, that is, without Theorem 3.3. Let us also mention some properties of the NEF $\mathcal{F}_{(\alpha,\beta),\mathrm{I}}$. From Theorem 2.1(1), we deduce that $L_{\mu_{(\alpha,\beta),G}}(y)$ is finite if and only if $-y \in P_G$. Since $-P_G$ is an open subset of $Z_G$ (see [12]), $\mathcal{F}_{(\alpha,\beta),\mathrm{I}}$ is a regular family in the sense of [4] and the domain of the means $M_{\mathcal{F}_{(\alpha,\beta),\mathrm{I}}}$ of the family $\mathcal{F}_{(\alpha,\beta),\mathrm{I}}$ coincides with the interior of the closed convex support of the family. Thus $M_{\mathcal{F}_{(\alpha,\beta),\mathrm{I}}} = Q_G$. The cumulant function of $\mu_{(\alpha,\beta),G}$ is

$$(4.20) \quad k_{\mu_{(\alpha,\beta),G}}(y) = \sum_{j=1}^{k} \alpha_j \log\det((-y^{-1})_{C_j}) - \sum_{j=2}^{k} \beta_j \log\det((-y^{-1})_{S_j}).$$

The computation of its differential requires some care and it is done in the Appendix in Propositions A.1 and A.2. We give the result here. If $X \sim W_{Q_G}(\alpha,\beta,\sigma)$, for $y = \hat{\sigma}^{-1} \in P_G$,

$$\begin{aligned}
&E_{W_{Q_G}(\alpha,\beta,\sigma)}(X) \\
(4.21) \quad &= \frac{d}{dy} k_{\mu_{(\alpha,\beta),G}}(y) \\
&= \sum_{j=1}^{k} \alpha_j \tau(\hat{\sigma}_{V\setminus C_j, C_j} \sigma_{C_j}^{-1}) \sigma_{C_j} - \sum_{j=2}^{k} \beta_j \tau(\hat{\sigma}_{V\setminus S_j, S_j} \sigma_{S_j}^{-1}) \sigma_{S_j},
\end{aligned}$$

where

$$\tau(\hat{\sigma}_{V\setminus C_j, C_j} \sigma_{C_j}^{-1}) \sigma_{C_j}$$

$$= \begin{pmatrix} I_{C_j} & 0 \\ \hat{\sigma}_{V\setminus C_j, C_j} \sigma_{C_j}^{-1} & I_{V\setminus C_j} \end{pmatrix} \begin{pmatrix} \sigma_{C_j} & 0 \\ 0 & 0 \end{pmatrix} \begin{pmatrix} I_{C_j} & \sigma_{C_j}^{-1} \tau \sigma_{C_j, V\setminus C_j} \\ 0 & I_{V\setminus C_j} \end{pmatrix}$$

(4.22)
$$= \begin{pmatrix} \sigma_{C_j} & \hat{\sigma}_{C_j,V\setminus C_j} \\ \hat{\sigma}_{V\setminus C_j,C_j} & \hat{\sigma}_{V\setminus C_j}\sigma_{C_j}^{-1}\sigma_{C_j,V\setminus C_j} \end{pmatrix}$$
$$= \hat{\sigma} - \hat{\sigma}_{V\setminus C_j \cdot C_j}$$

with a similar expression for $\tau(\sigma_{V\setminus S_j,S_j}\sigma_{S_j}^{-1})\sigma_{S_j}$. This implies that if $X \sim W_{Q_G}(\alpha,\beta,\sigma)$, then

$$E_{W_{Q_G}(\alpha,\beta,\sigma)}(X) = \left(\sum_{j=1}^k \alpha_j - \sum_{j=2}^k \beta_j\right)\hat{\sigma} + \sum_{j=1}^k \alpha_j \hat{\sigma}_{V\setminus C_j \cdot C_j} - \sum_{j=2}^k \beta_j \hat{\sigma}_{V\setminus S_j \cdot S_j}.$$

We have parallel results for the Type II Wishart. We note that for $(\alpha,\beta) \in \mathcal{B}$

$$\mathcal{F}_{(\alpha,\beta),\text{II}} = \{W_{P_G}(\alpha,\beta,\theta;dy), \theta \in Q_G\}$$

is the natural exponential family generated by the measure

$$\nu_{(\alpha,\beta),G}(dy) = \frac{H_G(\alpha,\beta;\varphi(y))}{\Gamma_{\text{II}}(\alpha,\beta)}\nu_G(dy).$$

Here the Laplace transform $L_{\nu_{(\alpha,\beta),G}}(x)$ is finite if and only if $x$ is in $-Q_G$ (Theorem 2.1, part 1.), which is an open subset of $I_G$ (see [12]), and the domain of the means of $\mathcal{F}_{\alpha,\beta,\text{II}}$ is $P_G$. We have $L_{\nu_{(\alpha,\beta),G}}(x) = H_G(\alpha,\beta;-x)$. This implies that the Laplace transform of $W_{P_G}(\alpha,\beta;\theta)$ is defined for $-x + \theta \in Q_G$ by

(4.23)
$$\int_{P_G} e^{\langle x,y\rangle} W_{P_G}(\alpha,\beta,\theta;dy) = \frac{H_G(\alpha,\beta;\theta-x)}{H_G(\alpha,\beta;\theta)}.$$

The cumulant transform is

$$k_{\nu_{(\alpha,\beta),G}}(x) = \sum_{j=1}^k \alpha_j \log\det((-x)_{C_j}) - \sum_{j=2}^k \beta_j \log\det((-x)_{S_j})$$

and its differential is given by the following element of $P_G$:

(4.24)
$$k'_\nu(x) = \sum_{j=1}^k \alpha_j (x_{C_j}^{-1})^0 - \sum_{j=2}^k \beta_j (x_{S_j}^{-1})^0.$$

This implies that if $Y \sim W_{P_G}(\alpha,\beta,\theta)$, then

(4.25)
$$E_{W_{P_G}(\alpha,\beta,\theta)}(Y) = \sum_{j=1}^k \alpha_j(-\theta_{C_j}^{-1})^0 - \sum_{j=2}^k \beta_j(-\theta_{S_j}^{-1})^0.$$

To conclude this section, let us make a few remarks. Formula (4.25) with $Y$ replaced by $\hat{\Sigma}_G^{-1}$, $(\alpha,\beta)$ replaced by $(\alpha'-\alpha,\beta'-\beta)$ and $\theta$ replaced by $D + x$ gives the posterior mean $E(\widehat{\Sigma}_G^{-1}|X = x)$ of the inverse of the natural parameter $\Sigma_G$ for $W_{Q_G}(\alpha, \beta, \Sigma_G)$ in Theorem 4.1 when the prior distribution on $\Sigma_G$ is $IW_{P_G}(\alpha', \beta', D)$. It is, of course, also of interest to compute the posterior mean of the natural parameter $\Sigma_G$. In other words, we need $E_{W_{P_G}(\alpha'-\alpha,\beta'-\beta,D+x)}(\varphi(Y))$ when $Y = \widehat{\Sigma}_G^{-1}$. We have the general formula

$$
\begin{aligned}
-\theta = &\left(\sum_{j=1}^{k}\left(\alpha_j + \frac{c_j+1}{2}\right) - \sum_{j=2}^{k}\left(\beta_j + \frac{s_j+1}{2}\right)\right) E_{W_{P_G}(\alpha,\beta,\theta)}(\varphi(Y)) \\
&+ \sum_{j=1}^{k}\left(\alpha_j + \frac{c_j+1}{2}\right) E_{W_{P_G}(\alpha,\beta,\theta)}(\varphi(Y)_{V \setminus C_j \cdot C_j}) \\
&- \sum_{j=2}^{k}\left(\beta_j + \frac{s_j+1}{2}\right) E_{W_{P_G}(\alpha,\beta,\theta)}(\varphi(Y)_{V \setminus S_j \cdot S_j}).
\end{aligned}
$$
(4.26)

The proof is not straightforward. To derive (4.26), we use Stokes' formula and obtain

$$0 = \int_{P_G} (u'(y)v(y) + u(y)v'(y))\,dy$$

with $u(y) = H_G(\alpha + \frac{c+1}{2}, \beta + \frac{s+1}{2}, \varphi(y))$, $v(y) = e^{\langle \theta, y \rangle}$, and where $(\alpha, \beta) \in \mathcal{B}$ satisfy some restrictions similar to the restrictions for the existence of the expected value of the inverse Wishart, that is $p > \frac{r+1}{2}$ for the $w_r(p, \sigma)$ distribution. Using the same substitutions for $\alpha, \beta, \theta$ in (4.26) as we did in (4.25), we see that (4.26) implies that in Theorem 4.1, the posterior mean of $\Sigma_G$ is not linear in $x$; this is in accordance with Theorem 3 of [9]. We also see from Corollary 4.1 that when $X = \pi(S) \in Q_G$ follows the hyper Wishart distribution with $(\alpha, \beta) = (\frac{n}{2}, \frac{n}{2})$ and the shape parameters of the prior on $\Sigma_G$ are $(\alpha', \beta')$, then the shape parameters of the posterior are $(\alpha' - \frac{n}{2}, \beta' - \frac{n}{2})$. That is, as for the inverse Wishart or the hyper inverse Wishart, the parameters $(\alpha', \beta')$ of the $IW_{P_G}$ are added to half of the sample size and from (4.26) the choice of $(\alpha', \beta')$ has the same kind of impact on the posterior mean as the choice of the shape parameters for the inverse or hyper inverse Wishart.

Let us also mention here that the $IW_{P_G}^{**}(\alpha, \beta, \theta)$ distribution as given by equations (4.10)–(4.12) is conditionally $(k+1)$-reducible and is an enriched conjugate family of prior distributions, in the sense of [5], for the parameter $\Sigma_G = \pi(\Sigma)$ of a Gaussian distribution Markov with respect to $G$. This follows immediately from Theorem 4.4. The $IW_{P_G}^{**}(\alpha, \beta, \theta)$ is also closely linked to the enriched standard conjugate Wishart family of priors for $K = \Sigma^{-1}$ in the standard Gaussian distribution, that is, when $G$ is a complete graph, built by Consonni and Veronese [6]. Theorem 2 and Corollary 1 in that paper correspond to Theorem 4.4 here. However, we should note that it is the $W_{P_G}(\alpha,\beta,\theta)$ family that is an exponential family, not the $IW^{**}_{P_G}(\alpha,\beta,\theta)$ family and therefore the analog of the enriched Wishart, for $G$ decomposable, is $W_{P_G}(\alpha,\beta,\theta)$.

**5. Open problems.** We will now raise some natural questions related to the paper.

*Singularity.* The well-known Gyndikin theorem states that the mapping $\theta \mapsto (\det(-\theta))^{-p}$ from $-M_r^+$ to $(0,\infty)$ is the Laplace transform of some positive measure $\mu_p$ on symmetric real matrices of order $n$ if and only if $p$ is in the set

$$\Lambda = \left\{\frac{1}{2}, \frac{2}{2}, \ldots, \frac{r-1}{2}\right\} \cup \left(\frac{r-1}{2}, \infty\right).$$

A very readable proof of this theorem can be found in [17]. The natural exponential family generated by $\mu_p$ is the set of Wishart distributions with shape parameter $p$. If $p = j/2$ with $j = 1, \ldots, n-1$, then $\mu_p$ is concentrated on the singular semipositive definite matrices of rank $j$. For a decomposable graph $G$ on $V = \{1, \ldots, r\}$ and for $(\alpha, \beta) \in \mathcal{A}$, the mapping $y \mapsto H_G(\alpha, \beta; -\varphi(y))$ from $-P_G$ to $(0, \infty)$ is the Laplace transform of a positive measure on $Q_G$ which generates the natural exponential family of Wishart distributions of Type I. Natural questions are:

- For which values of $\alpha, \beta$ is $y \mapsto H_G(\alpha, \beta; -\varphi(y))$ the Laplace transform of some positive measure on $I_G$?
- How do we describe these measures?

Similar questions arise with the Wishart distributions of Type II: for which values of $\alpha, \beta$ is the mapping $x \mapsto H_G(\alpha, \beta; -x)$ from $-Q_G$ to $(0, \infty)$ the Laplace transform of some positive measure on $Z_G$?

*Complex and quaternionic numbers.* Wishart matrices with complex and quaternionic entries are well defined. Thus many concepts of the present paper are extendable to complex or quaternionic matrices in a rather mechanical way.

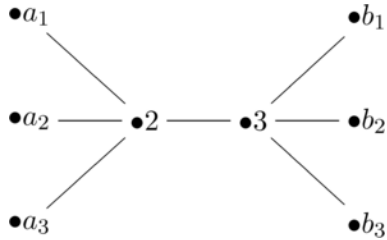

FIG. 3. *The case $n = m = 3$.*

*The sets $\mathcal{A}$ and $\mathcal{B}$.* Is it true that $\mathcal{A} = \bigcup A_P$ and that $\mathcal{B} = \bigcup B_P$ for any nonhomogeneous graph? Calculations are terrifying for the graph $A_5 : \bullet - \bullet - \bullet - \bullet - \bullet$. On the other hand, this conjecture is easily proved for the tree represented in Figure 3, with $n + m + 2$ vertices denoted $a_1, \ldots, a_n, b_1, \ldots, b_m, 2, 3$ with edges $a_i \sim 2$, $b_j \sim 3$ for all $i, j$ and $2 \sim 3$. To prove it, we need only extend the calculations of Proposition 3.2 and Corollary 3.1.

## APPENDIX

### A.1. Proofs of Section 2.

PROOF OF THEOREM 2.1(1). The dual of $P_G$ is $P_G^* = \{x \in I_G; \operatorname{tr}(xy) > 0$ for all $y \in \overline{P}_G \setminus \{0\}\}$ where $\overline{P}_G$ is the closure of $P_G$, that is, the cone of positive semidefinite matrices of $Z_G$. To show that $Q_G = P_G^*$ we will first show that $Q_G \subset P_G^*$ and then that $Q_G \supset P_G^*$. If $x \in Q_G$, then by Theorem 2.1 there exists a symmetric positive definite matrix $\hat{x}$ which is the completion of $x$. Thus $\operatorname{tr}(xy) = \operatorname{tr}(\hat{x}y)$ for all $y \in Z_G$. Furthermore, if $y \in \overline{P}_G \setminus \{0\}$ then $\operatorname{tr}(\hat{x}y) = \operatorname{tr}((\hat{x})^{1/2}y(\hat{x})^{1/2})$. Since the matrix $((\hat{x})^{1/2}y(\hat{x})^{1/2})$ is positive semidefinite and nonzero, its trace is positive. Thus $x \in P_G^*$ and $Q_G \subset P_G^*$ is proved.

Conversely, take $x \in I_G$ such that $x \in P_G^*$. Fix a clique $C$ and consider a vector $v$ of $\mathbb{R}^r$ such that the components of $v$ which are not in $C$ are 0. Denote by $v_C$ and by $x_C$ the restrictions of $v$ to $C$ and to $C \times C$, respectively, and assume that $v \neq 0$ and thus $v_C \neq 0$. Since $vv^t$ is in $\overline{P}_G \setminus \{0\}$ and since $x \in P_G^*$, $0 < \operatorname{tr}(xvv^t) = v^t xv = v_C^t x_C v_C$. Moreover, this is true for any $v_C \neq 0$ and therefore $x_C$ is positive definite. Since this is true for all cliques, we deduce that $x$ is in $Q_G$ and $Q_G \supset P_G^*$ is also proved. We thank S. Andersson for this result and this proof. Our former proof was longer and relied on the description of the extremal lines of the cones $\overline{P_G}$ and $\overline{Q_G}$ given in [12]. □

### A.2. Proofs of Section 3.3.

PROOF OF PROPOSITION 3.1. For convenience, we agree to write $t$ for a vertex of $T$ corresponding to a clique $C$ while we write $q$ for a vertex of $T$ corresponding to a separator $S$. As defined in (3.10), $m_v = \sum_{u \preceq v} n_u$ and for any $v \in T$, we have

$$|x_C| = \prod_{u \preceq t} |x_{[u]}|, \qquad |x_S| = \prod_{u \preceq q} |x_{[u]}|, \qquad |C| = c = m_t, \qquad |S| = s = m_q.$$

This means that $|x_{[u]}|$ will appear in $|x_C|$ for any $C$ such that $u \preceq t$ and in any $|x_S|$ for any $S$ such that $u \preceq q$. Therefore

$$\frac{\prod |x_C|^{\alpha_t - (c+1)/2}}{\prod |x_S|^{\nu(q)(\beta_q - (s+1)/2)}}$$

$$= \prod_{u \in T} |x_{[u]\cdot}|^{\sum_{u \preceq t}(\alpha_t - (m_t+1)/2) - \sum_{u \preceq q} \nu(q)(\beta_q - (m_q+1)/2)}$$

$$(A.1) = \prod_{u \in T} |x_{[u]\cdot}|^{(\sum_{u \preceq t} \alpha_t - \sum_{u \preceq q} \nu(q)\beta_q)} |x_{[u]\cdot}|^{-(\sum_{u \preceq t}(m_t+1)/2 - \sum_{u \preceq q} \nu(q)(m_q+1)/2)}$$

$$= \prod_{u \in T} |x_{[u]\cdot}|^{(\sum_{u \preceq t} \alpha_t - \sum_{u \preceq q} \nu(q)\beta_q)} |x_{[u]\cdot}|^{-(\sum_{u \prec t} n_t/2 + (m_u+1)/2)}$$

$$= \prod_{u \in T} |x_{[u]\cdot}|^{(\sum_{u \preceq t} \alpha_t - \sum_{u \preceq q} \nu(q)\beta_q)} |x_{[u]\cdot}|^{-(\sum_{u \prec v} n_v/2 + \sum_{v \prec u} n_v/2 + (n_u+1)/2)},$$

where the third equality above follows from the definition of $m_v$ and the fact that for any vertex $q$ of the tree, $\nu(q)$ is equal to the number of children of $q$ minus 1 (see part 2 of Proposition 2.2). We now make the change of variables (3.14). The Jacobian of this change of variables is

$$J = \prod_{v \in T} |x_{<v>}|^{n_v} = \prod_{v \in T} \left( \prod_{u \prec v} |x_{[u]\cdot}| \right)^{n_v}$$

$$(A.2)$$
$$= \prod_{u \in T} |x_{[u]\cdot}|^{\sum_{u \prec v} n_v}.$$

Therefore we obtain the image of $H_G(\alpha, \beta, x) \mu_G(dx)$ as

$$H_G^*(\alpha, \beta, x) \mu_G^* \left( \prod_{u \in T} dx_{[u]\cdot} \, d(x_{[u>} x_{<u>}^{-1}) \right)$$

$$= \prod_{t \in T} |x_{[t]\cdot}|^{\alpha_t - (m_t+1)/2}$$

$$\times \prod_{u \in T, u \neq t} |x_{[u]\cdot}|^{(\sum_{u \preceq t} \alpha_t - \sum_{u \preceq q} \nu(q)\beta_q + \sum_{u \prec v} n_v/2 - \sum_{v \prec u} n_v/2) - (n_u+1)/2}$$

$$(A.3) \qquad \times dx_{[u]\cdot} \, d(x_{[u>} x_{<u>}^{-1})$$

$$= \prod_{u \in T} |x_{[u]\cdot}|^{(\sum_{u \preceq t} \alpha_t - \sum_{u \preceq q} \nu(q)\beta_q + \sum_{u \prec v} n_v/2 - \sum_{v \prec u} n_v/2) - (n_u+1)/2}$$

$$\times dx_{[u]\cdot} \, d(x_{[u>} x_{<u>}^{-1})$$

$$= \prod_{u \in T} |x_{[u]\cdot}|^{\lambda_u - (n_u+1)/2} dx_{[u]\cdot} \, d(x_{[u>} x_{<u>}^{-1}),$$

where

$$\lambda_u = \sum_{u \preceq t} \alpha_t - \sum_{u \preceq q} \nu(q)\beta_q + \sum_{u \prec v} \frac{n_v}{2} - \sum_{v \prec u} \frac{n_v}{2}.$$

□

**A.3. Proofs of Section 3.4.** In the sequel, in order to avoid numbering difficulties for separators with multiplicity greater than one, we sometimes use the generic notation $S$ for a separator and $\nu(S)$ for its multiplicity. However, when it is important to list the separators as they appear from a perfect order of the cliques, we denote the separator $S_j$ by $\langle j \rangle$. The double notation should not cause any difficulty.

PROOF OF THEOREM 3.3. To avoid any ambiguity in the notation in the proof below and all the other proofs in the reminder of the paper, let us recall that the set $\mathcal{S}$ of distinct separators contains $k' \leq k-1$ elements. For convenience, let us write $A$ for the left-hand side of (3.19) and $y = \hat{\sigma}^{-1}$ for $\sigma \in Q_G$. Using the Jacobian (2.13) and formula (2.9), we then have

$$A = \int_{Q_G} e^{-\langle x, \hat{\sigma}^{-1}\rangle} \frac{\prod_{j=1}^{k} |x_{C_j}|^{\alpha_j - (c_j+1)/2}}{\prod_{S \in \mathcal{S}} |x_S|^{\nu(S)(\beta(S)-(|S|+1)/2)}} \, dx$$

$$= \int |x_{C_1}|^{\alpha_1 - (c_1+1)/2} e^{-\langle x_{C_1}, \sigma_{C_1}^{-1}\rangle}$$

$$\times \prod_{j=2}^{k} |x_{[j].}|^{\alpha_j - (c_j+1)/2} e^{-\langle x_{[j].}, \sigma_{[j].}^{-1}\rangle}$$

$$\times \prod_{j=2}^{k} e^{-\langle (x_{[j>}x_{<j>}^{-1} - \sigma_{[j>}\sigma_{<j>}^{-1}), \sigma_{[j].}^{-1}(x_{[j>}x_{<j>}^{-1} - \sigma_{[j>}\sigma_{<j>}^{-1})x_{<j>}\rangle}$$

$$\times \prod_{S \in \mathcal{S}} |x_S|^{\sum_{i \in J(P,S)}(\alpha_i - (c_i+1)/2) - \nu(S)(\beta(S)-(|S|+1)/2)}$$

$$\times \prod_{S \in \mathcal{S}} |x_S|^{\sum_{i \in J(P,S)} c_i - \nu(S)|S|} \, dx_{C_1} \prod_{j=2}^{k} d(x_{[j>}x_{<j>}^{-1}) \, dx_{[j].}.$$

Now, since the cardinality of $J(P,S)$ is equal to $\nu(S)$ and, by assumption 1 of the theorem, all $|x_{S_j}| = |x_{<j>}|, j \neq 2$, appear with exponent equal to $\frac{c_j - s_j}{2}$ while $|x_{S_2}| = |x_{<2>}|$ appears with exponent equal to $\sum_{i \in J(P,S_2)} \alpha_i - \nu(q)\beta_q + \sum_{i \in J(P,S_2)} \frac{c_i - s_2}{2}$, we have

$$A = \int |x_{C_1}|^{\alpha_1 - (c_1+1)/2} e^{-\langle x_{C_1}, \sigma_{C_1}^{-1}\rangle} |x_{<2>}|^{\delta_2}$$

(A.4) $\quad \times \prod_{j=2}^{k} |x_{<j>}|^{(c_j - s_j)/2} e^{-\langle (x_{[j>}x_{<j>}^{-1} - \sigma_{[j>}\sigma_{<j>}^{-1}), \sigma_{[j].}^{-1}(x_{[j>}x_{<j>}^{-1} - \sigma_{[j>}\sigma_{<j>}^{-1})x_{<j>}\rangle}$

$$\times |x_{[j].}|^{\alpha_j - s_j/2 - (c_j - s_j+1)/2} e^{-\langle x_{[j].}, \sigma_{[j].}^{-1}\rangle} \, dx_{C_1} \prod_{j=2}^{k} d(x_{[j>}x_{<j>}^{-1}) \, dx_{[j].}.$$

Clearly $(x_{[k].}, x_{[k>}x_{<k>}^{-1})$ is independent of $x_{C_1\cup\cdots\cup C_{k-1}} = (x_{C_1}, x_{[j>}x_{<j>}^{-1}, x_{[j].}, j = 2, \ldots, k-1)$ and $x_{[k].}$ is independent of $x_{[k>}x_{<k>}^{-1}$. Therefore holding all other variables fixed, we first integrate with respect to $(x_{[k].}, x_{[k>}x_{<k>}^{-1})$. Since $<k> \subset C_1 \cup \cdots \cup C_{k-1}$, then $x_{<k>}$ is fixed and by Lemma 2.4 we obtain

$$\int_{M^+_{c_k-s_k}} |x_{[k].}|^{\alpha_k - s_k/2 - (c_k - s_k + 1)/2} e^{-\langle x_{[k].}, \sigma^{-1}_{[k].}\rangle} \, dx_{[k].}$$

$$\times \int_{L(R^{s_k}, R^{c_k-s_k})} |x_{<k>}|^{(c_k-s_k)/2}$$

$$\times e^{-\langle (x_{[k>}x_{<k>}^{-1} - \sigma_{[k>}\sigma_{<k>}^{-1}), \sigma^{-1}_{[k].}(x_{[k>}x_{<k>}^{-1} - \sigma_{[k>}\sigma_{<k>}^{-1})x_{<k>}\rangle}$$

$$\times d(x_{[k>}x_{<k>}^{-1})$$

$$= \Gamma_{c_k - s_k}\left(\alpha_k - \frac{s_k}{2}\right)|\sigma_{[k].}|^{\alpha_k}$$

$$\times \int_{L(R^{s_k}, R^{c_k-s_k})} |\sigma_{[k].}|^{-s_k/2}|x_{<k>}|^{(c_k-s_k)/2}$$

$$\times e^{-\langle (x_{[k>}x_{<k>}^{-1} - \sigma_{[k>}\sigma_{<k>}^{-1}), \sigma^{-1}_{[k].}(x_{[k>}x_{<k>}^{-1} - \sigma_{[k>}\sigma_{<k>}^{-1})x_{<k>}\rangle}$$

$$\times d(x_{[k>}x_{<k>}^{-1})$$

$$= \pi^{(s_k(c_k-s_k))/2}\Gamma_{c_k-s_k}\left(\alpha_k - \frac{s_k}{2}\right)|\sigma_{[k].}|^{\alpha_k}.$$

Repeating this process successively for $j = k-1, \ldots, 2$, we obtain

$$A = \prod_{j=2}^{k} \pi^{s_j(c_j-s_j)/2}\Gamma_{c_j-s_j}\left(\alpha_j - \frac{s_j}{2}\right)|\sigma_{[j].}|^{\alpha_j}$$

$$\times \int |x_{C_1}|^{\alpha_1 - (c_1+1)/2} e^{-\langle x_{C_1}, \sigma^{-1}_{C_1}\rangle}|x_{<2>}|^{\delta_2} \, dx_{C_1}.$$

In this last integral, setting, as in Section 4.2,

$$x_{[1].} = x_{C_1 \setminus S_2} - x_{C_1 \setminus S_2, S_2} x_{S_2}^{-1} x_{S_2, C_1 \setminus S_2}, \qquad x_{[12>} = x_{C_1 \setminus S_2, S_2} x_{S_2}^{-1},$$

we make the change of variable

(A.5) $$x_{C_1} \mapsto (x_{[1].}, x_{[12>}, x_{<2>})$$

with Jacobian equal to $|x_{<2>}|^{c_1-s_2}$. Then, by (2.11)

$$\int |x_{C_1}|^{\alpha_1 - (c_1+1)/2} e^{-\langle x_{C_1}, \sigma^{-1}_{C_1}\rangle}|x_{<2>}|^{\delta_2} \, dx_{C_1}$$

(A.6) $$= \int |x_{[1].}|^{\alpha_1 - s_2/2 - (c_1 - s_2 + 1)/2} e^{-\langle x_{[1].}, \sigma^{-1}_{[1].}\rangle} \, dx_{[1].}.$$

$$\times \int_{M_{s_2}^+} |x_{<2>}|^{\alpha_1 + (c_1 - s_2)/2 + \delta_2 - (s_2+1)/2} e^{-\langle x_{<2>}, \sigma_{S_2}^{-1}\rangle}$$

$$\times \left( \int_{L(R^{s_2}, R^{c_1-s_2})} e^{-\langle (x_{[12>} - \sigma_{[12>}), \sigma_{[1]}^{-1}.(x_{[12>} - \sigma_{[12>})x_{<2>}\rangle} \, dx_{[12>} \right) dx_{<2>}$$

$$= \Gamma_{c_1-s_2}\left(\alpha_1 - \frac{s_2}{2}\right) |\sigma_{[1]}|^{\alpha_1} \pi^{s_2(c_1-s_2)/2}$$

$$\times \int_{M_{s_2}^+} |x_{<2>}|^{\alpha_1 - s_2/2 - (c_1-s_2+1)/2 + c_1 - s_2 - (c_1-s_2)/2 + \delta_2}$$

$$\times e^{-\langle x_{<2>}, \sigma_{S_2}^{-1}\rangle} dx_{<2>}$$

$$= \pi^{s_2(c_1-s_2)/2} \Gamma_{c_1-s_2}\left(\alpha_1 - \frac{s_2}{2}\right) |\sigma_{[1]}|^{\alpha_1}$$

$$\times \int_{M_{s_2}^+} |x_{<2>}|^{\alpha_1 - (s_2+1)/2 + \delta_2} e^{-\langle x_{<2>}, \sigma_{S_2}^{-1}\rangle} dx_{<2>}$$

$$= \pi^{s_2(c_1-s_2)/2} \Gamma_{c_1-s_2}\left(\alpha_1 - \frac{s_2}{2}\right) |\sigma_{[1]}|^{\alpha_1} |\sigma_{S_2}|^{\alpha_1} |\sigma_{S_2}|^{\delta_2} \Gamma_{s_2}(\alpha_1 + \delta_2)$$

$$(A.7) \quad = \pi^{s_2(c_1-s_2)/2} |\sigma_{C_1}|^{\alpha_1} |\sigma_{S_2}|^{\delta_2} \Gamma_{c_1-s_2}\left(\alpha_1 - \frac{s_2}{2}\right) \Gamma_{s_2}(\alpha_1 + \delta_2).$$

Now, let us observe that using the multiplicity of $S_2$ and assumption 1 of the theorem, we obtain

$$(A.8) \quad \prod_{j=2}^{k} |\sigma_{[j]}|^{\alpha_j} = \frac{|\sigma_{C_2}|^{\alpha_2}}{|\sigma_{S_2}|^{\alpha_2}} \prod_{j=3}^{k} \frac{|\sigma_{C_j}|^{\alpha_j}}{|\sigma_{S_j}|^{\alpha_j}}$$

$$= \frac{|\sigma_{C_2}|^{\alpha_2}}{|\sigma_{S_2}|^{\sum_{i \in J(P, S_2)} \alpha_i}} \frac{\prod_{j=3}^{k} |\sigma_{C_j}|^{\alpha_j}}{\prod_{S \in \mathcal{S}, S \neq S_2} |\sigma_S|^{\nu(S)\beta(S)}}.$$

Combining (A.5), (A.7) and (A.8), we obtain (3.19) with

$$\Gamma_I(\alpha, \beta) = \Gamma_{s_2}(\alpha_1 + \delta_2)$$

$$\times \pi^{(1/2)(c_1 - s_2)s_2} \Gamma_{c_1 - s_2}\left(\alpha_1 - \frac{s_2}{2}\right)$$

$$\times \prod_{j=2}^{k} \pi^{\sum_{j=2}^{k}(1/2)(c_j - s_j)s_j} \Gamma_{c_j - s_j}\left(\alpha_j - \frac{s_j}{2}\right).$$

To obtain (3.20), we use (2.16). □

PROOF OF THEOREM 3.4. For convenience let us denote by $B$ the left-hand side of (3.22). Using first $\varphi(y) = x \in Q_G$ and the Jacobian (2.4) and then, making the change of variable (2.12) with Jacobian (2.13) and using (2.9), we have

$$B = \int_{Q_G} e^{-\langle \theta, \hat{x}^{-1} \rangle} \frac{\prod_{j=1}^{k} |x_{C_j}|^{\alpha_j - (c_j+1)/2}}{\prod_{S \in \mathcal{S}} |x_S|^{\nu(S)(\beta(S) - (|S|+1)/2)}} \, dx$$

$$= \int |x_{C_1}|^{\alpha_1 - (c_1+1)/2} e^{-\langle x_{C_1}^{-1}, \theta_{C_1} \rangle}$$

$$\times \prod_{j=2}^{k} |x_{[j]\cdot}|^{\alpha_j - (c_j+1)/2} e^{-\langle x_{[j]\cdot}^{-1}, \theta_{[j]\cdot} \rangle}$$

$$\times \prod_{j=2}^{k} e^{-\langle (x_{[j>}x_{<j>}^{-1} - \theta_{[j>}\theta_{<j>}^{-1}), x_{[j]\cdot}^{-1}(x_{[j>}x_{<j>}^{-1} - \theta_{[j>}\theta_{<j>}^{-1})\theta_{<j>} \rangle}$$

$$\times \prod_{S \in \mathcal{S}} |x_S|^{\sum_{i \in J(P,S)} (\alpha_i - (c_i+1)/2) - \nu(S)(\beta(S) - (|S|+1)/2)}$$

$$\times \prod_{S \in \mathcal{S}} |x_S|^{\sum_{i \in J(P,S)} c_i - \nu(S)|S|} \, dx_{C_1}$$

$$\times \prod_{j=2}^{k} d(x_{[j>}x_{<j>}^{-1}) \, dx_{[j]\cdot}.$$

By assumption 1 of the theorem, this is equal to

(A.9)
$$B = \int |x_{C_1}|^{\alpha_1 - (c_1+1)/2} e^{-\langle x_{C_1}^{-1}, \theta_{C_1} \rangle} |x_{<2>}|^{\gamma_2}$$

$$\times \prod_{j=2}^{k} |x_{[j]\cdot}|^{\alpha_j - (c_j+1)/2} e^{-\langle x_{[j]\cdot}^{-1}, \theta_{[j]\cdot} \rangle}$$

$$\times e^{-\langle (x_{[j>}x_{<j>}^{-1} - \theta_{[j>}\theta_{<j>}^{-1}), x_{[j]\cdot}^{-1}(x_{[j>}x_{<j>}^{-1} - \theta_{[j>}\theta_{<j>}^{-1})\theta_{<j>} \rangle} \, dx_{C_1}$$

$$\times \prod_{j=2}^{k} d(x_{[j>}x_{<j>}^{-1}) \, dx_{[j]\cdot}.$$

Clearly $x_{C_1}, (x_{[j>}x_{<j>}^{-1}, x_{[j]\cdot}, j = 2, \ldots, k)$ are mutually independent and

$$B = B_1 \times \prod_{j=2}^{k} B_j,$$

where

$$B_1 = \int |x_{C_1}|^{\alpha_1 - (c_1+1)/2} e^{-\langle x_{C_1}^{-1}, \theta_{C_1} \rangle} |x_{<2>}|^{\gamma_2} \, dx_{C_1}$$

and, using Lemma 2.1 in the third equality below, we have

$$B_j = \int_{M^+_{c_j-s_j}} \left( \int_{L(R^{s_j}, R^{c_j-s_j})} e^{-\langle (x_{[j>}x^{-1}_{<j>} - \theta_{[j>}\theta^{-1}_{<j>}), x^{-1}_{[j]\cdot}(x_{[j>}x^{-1}_{<j>} - \theta_{[j>}\theta^{-1}_{<j>})\theta_{<j>} \rangle} \right.$$

$$\left. \times d(x_{[j>}x^{-1}_{<j>}) \right)$$

$$\times |x_{[j]\cdot}|^{\alpha_j - (c_j+1)/2} e^{-\langle x^{-1}_{[j]\cdot}, \theta_{[j]\cdot} \rangle} dx_{[j]\cdot}.$$

$$= \pi^{s_j(c_j-s_j)/2} \int_{M^+_{c_j-s_j}} |\theta_{<j>}|^{-(c_j-s_j)/2} |x_{[j]\cdot}|^{\alpha_j - (c_j+1)/2 + s_j/2} e^{-\langle x^{-1}_{[j]\cdot}, \theta_{[j]\cdot} \rangle} dx_{[j]\cdot}.$$

$$= \pi^{s_j(c_j-s_j)/2} \int_{M^+_{c_j-s_j}} |\theta_{<j>}|^{-(c_j-s_j)/2} |x^{-1}_{[j]\cdot}|^{-\alpha_j + (c_j-s_j+1)/2 - (c_j-s_j+1)}$$

$$\times e^{-\langle x^{-1}_{[j]\cdot}, \theta_{[j]\cdot} \rangle} d(x^{-1}_{[j]\cdot})$$

$$= \pi^{s_j(c_j-s_j)/2} |\theta_{<j>}|^{-(c_j-s_j)/2} |\theta_{[j]\cdot}|^{\alpha_j} \Gamma_{c_j-s_j}(-\alpha_j).$$

Therefore

$$(A.10) \quad B = \prod_{j=2}^{k} \pi^{s_j(c_j-s_j)/2} |\theta_{<j>}|^{-(c_j-s_j)/2} |\theta_{[j]\cdot}|^{\alpha_j} \Gamma_{c_j-s_j}(-\alpha_j) \times B_1.$$

To compute $B_1$, let us make the change of variable (A.5). Then

$$B_1 = \int |x_{[1]\cdot}|^{\alpha_1 - (c_1+1)/2} e^{-\langle x^{-1}_{[1]\cdot}, \theta_{[1]\cdot} \rangle} e^{-\langle (x_{[12>} - \theta_{[12>}), x^{-1}_{[1]\cdot}(x_{[12>} - \theta_{[12>})\theta_{<2>} \rangle}$$

$$(A.11) \quad \times e^{-\langle x^{-1}_{<2>}, \theta_{<2>} \rangle} |x_{<2>}|^{\alpha_1 - (c_1-s_2)/2 - (s_2+1)/2 + \gamma_2 + (c_1-s_2)}$$

$$\times dx_{[1]\cdot} dx_{[12>} dx_{<2>}.$$

Integrating with respect to $x_{[12>}$ and using Lemma 2.1, we obtain

$$B_1 = |\theta_{<2>}|^{-(c_1-s_2)/2} \pi^{(c_1-s_2)s_2/2}$$

$$\times \int_{M^+_{c_1-s_2}} e^{-\langle x^{-1}_{[1]\cdot}, \theta_{[1]\cdot} \rangle} |x_{[1]\cdot}|^{\alpha_1 - (c_1-s_2+1)/2} dx_{[1]\cdot}$$

$$\times \int_{M^+_{s_2}} e^{-\langle x^{-1}_{<2>}, \theta_{<2>} \rangle}$$

$$\times |x_{<2>}|^{\alpha_1 - (c_1-s_2)/2 + \gamma_2 - (s_2+1)/2 + (c_1-s_2)} dx_{<2>}$$

$$= |\theta_{<2>}|^{-(c_1-s_2)/2} \pi^{(c_1-s_2)s_2/2}$$

$$\times \int_{M^+_{c_1-s_2}} e^{-\langle x^{-1}_{[1]\cdot}, \theta_{[1]\cdot} \rangle} |x^{-1}_{[1]\cdot}|^{-\alpha_1 - (c_1-s_2+1)/2} d(x^{-1}_{[1]\cdot})$$

(A.12)
$$\times \int_{M_{s_2}^+} e^{-\langle x_{<2>}^{-1}, \theta_{<2>}\rangle}$$
$$\times |x_{<2>}^{-1}|^{-\alpha_1-(c_1-s_2)/2-\gamma_2-(s_2+1)/2} d(x_{<2>}^{-1})$$
$$= \pi^{(c_1-s_2)s_2/2}|\theta_{<2>}|^{-(c_1-s_2)/2}|\theta_{[1]\cdot}|^{\alpha_1}$$
$$\times \Gamma_{c_1-s_2}(-\alpha_1)|\theta_{<2>}|^{\alpha_1+(c_1-s_2)/2+\gamma_2}\Gamma_{s_2}\left(-\alpha_1 - \frac{c_1-s_2}{2} - \gamma_2\right)$$
$$= \pi^{(c_1-s_2)s_2/2}|\theta_{<2>}|^{\alpha_1+\gamma_2}|\theta_{[1]\cdot}|^{\alpha_1}\Gamma_{c_1-s_2}(-\alpha_1)\Gamma_{s_2}\left(-\alpha_1 - \frac{c_1-s_2}{2} - \gamma_2\right)$$
$$= \pi^{(c_1-s_2)s_2/2}|\theta_{C_1}|^{\alpha_1}|\theta_{<2>}|^{\alpha_1}\Gamma_{c_1-s_2}(-\alpha_1)\Gamma_{s_2}\left(-\alpha_1 - \frac{c_1-s_2}{2} - \gamma_2\right).$$

Let us now observe that

$$\prod_{j=2}^{k}|\theta_{<j>}|^{-(c_j-s_j)/2}|\theta_{[j]\cdot}|^{\alpha_j}$$

(A.13)
$$= \frac{\prod_{j=2}^{k}|\theta_{C_j}|^{\alpha_j}|\theta_{<j>}|^{-\alpha_j-(c_j-s_j)/2}\prod_{S\in\mathcal{S}}|\theta_S|^{\nu(S)\beta(S)}}{\prod_{S\in\mathcal{S}}|\theta_S|^{\nu(S)\beta(S)}}$$

$$= \frac{\prod_{j=2}^{k}|\theta_{C_j}|^{\alpha_j}}{\prod_{S\in\mathcal{S}}|\theta_S|^{\nu(S)\beta(S)}}\prod_{S\in\mathcal{S}}|\theta_S|^{-\sum_{j\in J(P,S)}\alpha_j+(1/2)(c_j-s_j)+\nu(S)\beta(S)}$$

$$= \frac{\prod_{j=2}^{k}|\theta_{C_j}|^{\alpha_j}}{\prod_{S\in\mathcal{S}}|\theta_S|^{\nu(S)\beta(S)}}|\theta_{<2>}|^{-\gamma_2}.$$

Combining (A.10), (A.12) and (A.13), we obtain

$$\Gamma_{\mathrm{II}}(\alpha,\beta) = \pi^{(1/2)((c_1-s_2)s_2+\sum_{j=2}^{k}(c_j-s_j)s_j)}$$
$$\times \Gamma_{s_2}\left[-\alpha_1 - \frac{c_1-s_2}{2} - \gamma_2\right]\Gamma_{c_1-s_2}(-\alpha_1)\prod_{j=2}^{k}\Gamma_{c_j-s_j}(-\alpha_j).$$

To obtain (3.23), we use (2.17). □

### A.4. Proofs of Section 4.3.

PROPOSITION A.1. *For* $-y \in M_r^+$ *and* $C \subset \{1,\ldots,r\}$ *denote* $\sigma_C(y) = ((-y)^{-1})_C$. *We write* $y$ *by blocks corresponding to* $C$ *and its complement*

$$y = \begin{bmatrix} y_1 & y_{12} \\ y_{21} & y_2 \end{bmatrix}.$$

We denote for simplicity $y'_{12} = y_{12} y_2^{-1}$. With this notation the differential of $y \mapsto \sigma_C(y)$ is

$$(A.14) \quad h = \begin{bmatrix} h_1 & h_{12} \\ h_{21} & h_2 \end{bmatrix} \mapsto [\sigma_C \quad -\sigma_C y'_{12}] \begin{bmatrix} h_1 & h_{12} \\ h_{21} & h_2 \end{bmatrix} \begin{bmatrix} \sigma_C \\ -y'_{21} \sigma_C \end{bmatrix}.$$

Furthermore, the differential of $y \mapsto \kappa_C(y) = -\log \det \sigma_C(y)$ is

$$(A.15) \quad \begin{aligned} h &\mapsto \operatorname{tr} \begin{bmatrix} h_1 & h_{12} \\ h_{21} & h_2 \end{bmatrix} \begin{bmatrix} \sigma_C & -\sigma_C y'_{12} \\ -y'_{21} \sigma_C & y'_{21} \sigma_C y'_{12} \end{bmatrix} \\ &= \operatorname{tr} \begin{bmatrix} h_1 & h_{12} \\ h_{21} & h_2 \end{bmatrix} \begin{bmatrix} \sigma_C & \hat{\sigma}_{C, V \setminus C} \\ \hat{\sigma}_{V \setminus C, C} & \hat{\sigma}_{V \setminus C, C} \sigma_C^{-1} \hat{\sigma}_{C, V \setminus C} \end{bmatrix}, \end{aligned}$$

where the last equality is due to the fact that $y'_{12} = -\sigma_C^{-1} \hat{\sigma}_{C, V \setminus C}$.

PROOF. We know that

$$\sigma_C(y) = (y^{-1})_1 = -(y_1 - y_{12} y_2^{-1} y_{21})^{-1}.$$

Let $M_C$ and $M_C^+$ denote the restrictions of $M$ and $M_r^+$ to the clique $C$. Then $\sigma_C(y) = a \circ b(y)$ where $a: -M_C^+ \to M_C^+$ is defined by $a(x) = -x^{-1}$ and has differential $h \mapsto a'(x)(h) = x^{-1} h x^{-1}$ (a linear application from $M_C$ to $M_C$) and where $b: -M^+ \to -M_C^+$ is defined by $b(y) = y_1 - y_{12} y_2^{-1} y_{21}$. The differential of $b$ is the following linear mapping from $M$ to $M_C$:

$$(A.16) \quad \begin{aligned} \begin{bmatrix} h_1 & h_{12} \\ h_{21} & h_2 \end{bmatrix} &\mapsto h_1 - h_{12} y_2^{-1} y_{21} - y_{12} y_2^{-1} h_{21} + y_{12} y_2^{-1} h_2 y_2^{-1} y_{21} \\ &= [1 \quad -y'_{12}] \begin{bmatrix} h_1 & h_{12} \\ h_{21} & h_2 \end{bmatrix} \begin{bmatrix} 1 \\ -y'_{21} \end{bmatrix}. \end{aligned}$$

Finally we apply the composition of differentials to obtain

$$\sigma'_C(y)(h) = (a \circ b)'(y)(h) = a'(b(y))(b'(y)(h)) = \sigma_C(y)(b'(y)(h))\sigma_C(y),$$

which gives (A.14) when combined with (A.16). Now consider the real function $l$ defined on $M_C^+$ by $l(x) = \log \det x$. Then its differential is the linear form on $M_C$ defined by $h \mapsto l'(x)(h) = \operatorname{tr}(x^{-1} h)$. Thus the differential of the real function on $M_r^+$ defined by $l \circ \sigma_C$ is the following linear form on $M$:

$$h \mapsto (l \circ \sigma_C)'(y)(h) = \operatorname{tr} \sigma_C^{-1} \sigma'_C(y)(h),$$

which gives (A.15) when combined with (A.14). □

We will now use the previous proposition to compute $m = k'_\mu(y)$. We need to introduce the notation $h_C$ for the restriction of $h \in Z_G$ to $C \times C$, when $C \subset V$ the notation $C' = V \setminus C$ and the notation $h_{C,C'}$ and $h_{C',C}$ for the restrictions of $h$ to $C \times C'$ and $C' \times C$, respectively.

PROPOSITION A.2. *The differential of the real function $y \mapsto k_\mu(y)$ defined on $P_G$ by* (4.20) *is the linear form on $Z_G$ defined by*

$$
\begin{aligned}
h \mapsto & \sum_{j=1}^{k} \alpha_j [\operatorname{tr}(h_{C_j} \sigma_{C_j}) - 2\operatorname{tr}(h_{C'_j, C_j} \sigma_{C_j} y'_{C_j, C'_j}) + \operatorname{tr}(h_{C'_j} y'_{C'_j, C_j} \sigma_{C_j} y'_{C_j, C'_j})] \\
& - \sum_{j=2}^{k} \beta_j [\operatorname{tr}(h_{S_j} \sigma_{S_j}) - 2\operatorname{tr}(h_{S'_j, S_j} \sigma_{S_j} y'_{S_j, S'_j}) + \operatorname{tr}(h_{S'_j} y'_{S'_j, S_j} \sigma_{S_j} y'_{S_j, S'_j})] \\
= & \sum_{j=1}^{k} \alpha_j \operatorname{tr} \begin{pmatrix} h_{C_j} & h_{C_j, V \setminus C_j} \\ h_{V \setminus C_j, C_j} & h_{V \setminus C_j} \end{pmatrix} \begin{pmatrix} \sigma_{C_j} & \hat{\sigma}_{C_j, V \setminus C_j} \\ \hat{\sigma}_{V \setminus C_j, C_j} & \hat{\sigma}_{V \setminus C_j, C_j} \sigma_{C_j}^{-1} \hat{\sigma}_{C_j, V \setminus C_j} \end{pmatrix} \\
& - \sum_{j=2}^{k} \beta_j \operatorname{tr} \begin{pmatrix} h_{S_j} & h_{S_j, V \setminus S_j} \\ h_{V \setminus S_j, S_j} & h_{V \setminus S_j} \end{pmatrix} \begin{pmatrix} \sigma_{S_j} & \hat{\sigma}_{S_j, V \setminus S_j} \\ \hat{\sigma}_{V \setminus S_j, S_j} & \hat{\sigma}_{V \setminus S_j, S_j} \sigma_{S_j}^{-1} \hat{\sigma}_{S_j, V \setminus S_j} \end{pmatrix}.
\end{aligned}
$$

PROOF. We apply (A.15) to each term of the sum $k_\mu$. The proposition has been established for the cone $M_r^+$ and we apply it here to the restriction $P_G = M_r^+ \cap Z_G$ of $M_r^+$. Therefore the formulas for the differentials of functions restricted to this subspace are still in force when interpreted as linear applications defined on $Z_G$. □

**Acknowledgments.** The authors would like to thank S. A. Andersson for many lively discussions that led to significant improvements in this paper, S. L. Lauritzen for telling the authors that some hyper Markov property must hold for the Type I and II Wisharts, and finally a referee and the Associate Editor for reading the manuscript meticulously and suggesting several improvements in the presentation.

## REFERENCES

[1] ABRAMOWITZ, M. and STEGUN, I. A., eds. (1965). *Handbook of Mathematical Functions.* U.S. Government Printing Office, Washington.
[2] ANDERSSON, S. A. (1999). Private communication.
[3] ANDERSSON, S. A. and WOJNAR, G. G. (2004). Wishart distributions on homogeneous cones. *J. Theoret. Probab.* **17** 781–818. MR2105736
[4] BARNDORFF-NIELSEN, O. (1978). *Information and Exponential Families in Statistical Theory.* Wiley, Chichester. MR0489333
[5] CONSONNI, G. and VERONESE, P. (2001). Conditionally reducible natural exponential families and enriched conjugate priors. *Scand. J. Statist.* **28** 377–406. MR1842256
[6] CONSONNI, G. and VERONESE, P. (2003). Enriched conjugate and reference priors for the Wishart family on symmetric cones. *Ann. Statist.* **31** 1491–1516. MR2012823
[7] DAWID, A. and LAURITZEN, S. L. (1993). Hyper-Markov laws in the statistical analysis of decomposable graphical models. *Ann. Statist.* **21** 1272–1317. MR1241267
[8] DEMPSTER, A. P. (1972). Covariance selection. *Biometrics* **28** 157–175.

[9] DIACONIS, P. and YLVISAKER, D. (1979). Conjugate priors for exponential families. *Ann. Statist.* **7** 269–281. MR0520238

[10] GRÖNE, R., JOHNSON, C. R., DE SÁ, E. and WOLKOWICZ, H. (1984). Positive definite completions of partial Hermitian matrices. *Linear Algebra Appl.* **58** 109–124. MR0739282

[11] LAURITZEN, S. L. (1996). *Graphical Models.* Oxford Univ. Press. MR1419991

[12] LETAC, G. and MASSAM, H. (2006). Extremal rays and duals for cones of positive definite matrices with prescribed zeros. *Linear Algebra Appl.* **418** 737–750. MR2260225

[13] MASSAM, H. and NEHER, E. (1997). On transformations and determinants of Wishart variables on symmetric cones. *J. Theoret. Probab.* **10** 867–902. MR1481652

[14] MUIRHEAD, R. (1982). *Aspects of Multivariate Statistical Theory.* Wiley, New York. MR0652932

[15] OLKIN, I. and RUBIN, H. (1964). Multivariate beta distributions and independence properties of the Wishart distribution. *Ann. Math. Statist.* **35** 261–269. MR0160297

[16] ROVERATO, A. (2000). Cholesky decomposition of a hyper inverse Wishart matrix. *Biometrika* **87** 99–112. MR1766831

[17] SHANBAG, D. N. (1988). The Davidson-Kendall problem and related results on the structure of the Wishart distribution. *Austral. J. Statist.* **30A** 272–280.

[18] WHITTAKER, J. (1990). *Graphical Models in Applied Multivariate Statistics.* Wiley, Chichester. MR1112133

LABORATOIRE DE STATISTIQUE
ET PROBABILITÉS
UNIVERSITÉ PAUL SABATIER
31062 TOULOUSE
FRANCE
E-MAIL: letac@cict.fr

DEPARTMENT OF MATHEMATICS
AND STATISTICS
YORK UNIVERSITY
TORONTO, ONTARIO
CANADA M3J 1P3
E-MAIL: massamh@yorku.ca